\documentclass[sigconf,authorversion,nonacm]{acmart}

\usepackage{amsmath}
\usepackage{algorithmic}
\usepackage{graphicx}
\usepackage{textcomp}
\def\BibTeX{{\rm B\kern-.05em{\sc i\kern-.025em b}\kern-.08em
    T\kern-.1667em\lower.7ex\hbox{E}\kern-.125emX}}

\usepackage[linesnumbered,lined,ruled,commentsnumbered]{algorithm2e} 
\usepackage{subfigure}
\usepackage{multirow}

\let\oldnl\nl
\newcommand{\nonl}{\renewcommand{\nl}{\let\nl\oldnl}}

\DeclareMathOperator*{\argmin}{arg\,min}

\theoremstyle{definition}
\newtheorem{definition}{{Definition}}
\newcounter{examplecounter}
\newcommand{\example}[1]{\refstepcounter{examplecounter}\label{#1}{\textit{Example \arabic{examplecounter}:}}}
\newcounter{lemmacounter}
\newcommand{\lemma}[1]{\refstepcounter{lemmacounter}\label{#1}{\textsc{Lemma \arabic{lemmacounter}:}}}

\copyrightyear{}
\acmYear{}

\acmBooktitle{MobiQuitous 2020 - 17th EAI International Conference on Mobile and
Ubiquitous Systems: Computing, Networking and Services (MobiQuitous '20),
December 7--9, 2020, Darmstadt, Germany}
\acmPrice{15.00}
\acmDOI{10.1145/3448891.3448912}
\acmISBN{}

\begin{document}

\title{dFDA-VeD: A Dynamic Future Demand Aware Vehicle Dispatching System}

\author{Yang Guo$^1$, Tarique Anwar$^1$, Jian Yang$^{1, 2}$, and Jia Wu$^1$}
\affiliation{%
  \institution{1. Department of Computing, Macquarie University, Australia\\
  2. College of Computer Science and Technology, Dong Hua University, P.R. China.\\}
}
\email{{yang.guo5@students., tarique.anwar@, jian.yang@, jia.wu@}mq.edu.au}

\renewcommand{\shortauthors}{Guo et al.}


\begin{abstract}
With the rising demand of smart mobility, ride-hailing service is getting popular in the urban regions. These services maintain a system for serving the incoming trip requests by dispatching available vehicles to the pickup points. As the process should be socially and economically profitable, the task of vehicle dispatching is highly challenging, specially due to the time-varying travel demands and traffic conditions. Due to the uneven distribution of travel demands, many idle vehicles could be generated during the operation in different subareas. Most of the existing works on vehicle dispatching system, designed static relocation centers to relocate idle vehicles. However, as  traffic conditions and demand distribution dynamically change over time, the static solution can not fit the evolving situations. In this paper, we propose a dynamic future demand aware vehicle dispatching system. It can dynamically search the relocation centers considering both travel demand and traffic conditions. We evaluate the system on real-world dataset, and compare with the existing state-of-the-art methods in our experiments in terms of several standard evaluation metrics and operation time. Through our experiments, we demonstrate that the proposed system significantly improves the serving ratio and with a very small increase in operation cost.
\end{abstract}





\maketitle


\section{Introduction}

\textit{Mobility} and \textit{mobility-on-demand} services are major concerns of smart transportation. Mobility services exist in  cities since long in the form of public transportation, and mobility-on-demand services were limited to renting cars, offered by companies like Hertz and Avis. However, as we strive to make our cities smarter, over the past 10 years, we have seen the evolution of mobility-on-demand services into many new, effective, and more convenient forms. These services are increasingly being promoted as an influential strategy to address the challenges of urban transportation in large and fast-growing cities. Evolving from the traditional taxi service model, today  companies like Uber, Lyft, Ola, Didi and many others are popular as \textit{ride-hailing} mobility-on-demand service providers in many cities globally. These services are facilitated by the recent advancements in communication technologies and widely used GPS-enabled mobile devices. Customers can send trip requests to the service provider from their mobile devices in real time. Upon receiving the requests, the vehicle dispatching system of the service provider assigns available vehicles to serve the trip requests. While in progress, the system keeps track of geographic location of both the customer and the vehicle in order to maintain an updated information dynamically, and provide a smooth service. One major problem in the vehicle dispatching system is to find the most suitable vehicle to serve a trip request in such a way that results into the highest overall social and economical benefit. With rapid developments taking place currently in the field of Internet of Things (IoT), especially vehicle-to-anything (V2X), the future will have the availability of more data of vehicles and traffic, with high accuracy \cite{chen2017vehicle}. Using the dynamic traffic data on the roads and predicting the future travel demands are major aids in making estimations of benefit while serving a trip request, and thus have a high potential to improve the decision making of mobility-on-demand services.

\begin{figure}[!ht]
\vspace{-6pt}
\centering
\includegraphics[width=1\linewidth]{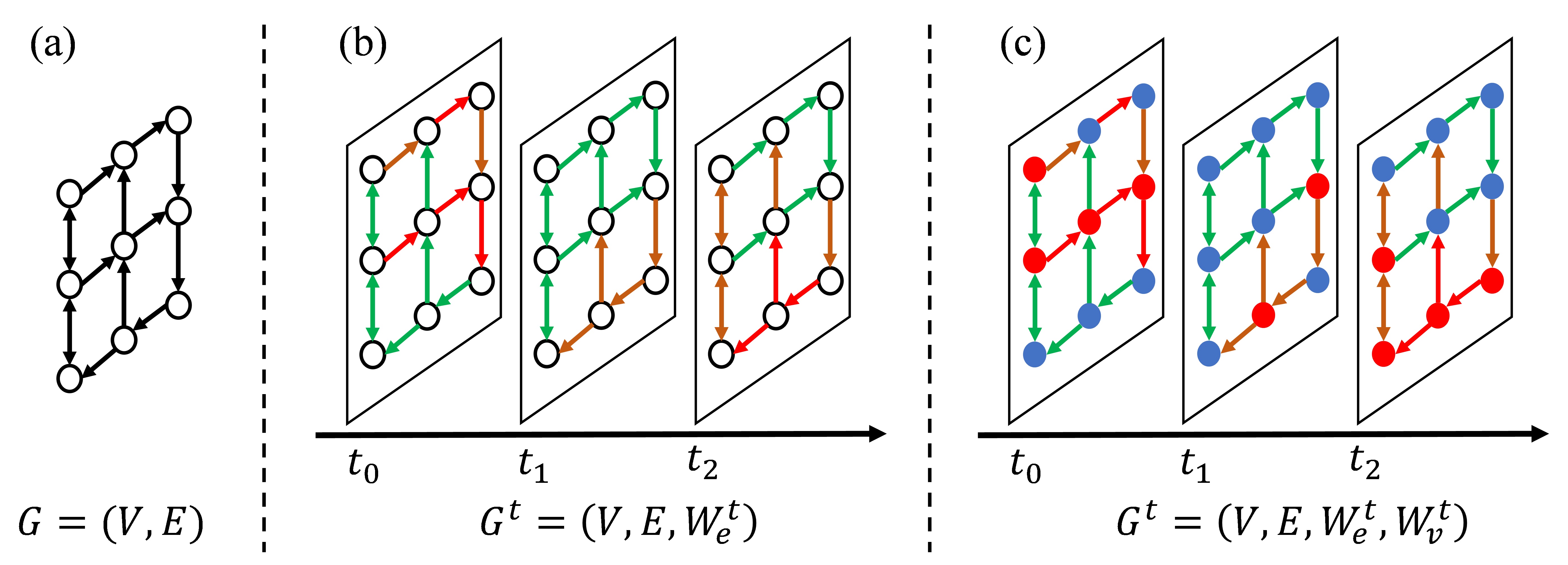}
\vspace{-8pt}
\caption{\textbf{Illustration of static and dynamic road graphs.} The directed graphs show the road network of a serving area, where vertices represent the pickup/dropoff points, and edges represent the directed path connecting the pickup/dropoff points. The colors highlight the different attributes on vertices and edges. The vertices in red color denote the points that have under-supplied vehicles and the vertices in blue color denote the points that over-supplied vehicles. The color of edges denotes the travel time in respective ways. (a) Static Road Graph $G= (V, E)$. The attributes of both vertices and edges do not change over time in this graph. (b) Dynamic-edge Road Graph $G^t = (V, E, W_e^t)$. The attributes of edges change over time in this graph. The changing color of edges in different timestamps show the time-varying attributes $W_e^t$ of edges. (c) Dynamic Graph $G^t= (V_t, E_t, W_v^t, W_e^t)$. The attributes of both vertices and edges change over time.}
\vspace{-6pt}
\label{fig:dynamicgraphs}
\end{figure}

Being an important component of the mobility-on-demand services, research on vehicle dispatching system has been under attention since quite some time \cite{dantzig1959truck, toth2002vehicle, pisinger2007general, kool2018attention}. Several concerns have been studied since then \cite{wang2019ridesourcing}. 
One major concern is to deal with the uneven geographic distribution of trip requests and available vehicles. Quite often there arise geographic areas of low trip demand with over-supplied vehicles, at the same time when there are other areas of high trip demand with under-supplied vehicles. To address this issue, efficient relocation of idle vehicles from areas of low demand to those of high demand has been an important problem \cite{guo2020FDAVeD}. Addressing this issue is crucial to improve the serving ratio of incoming trip requests. Vazifeh et al. \cite{vazifeh2018addressing} gave the lower bound of fleet size to serve all travel requests in an ideal scenario, considering that all travel requests are known in advance. Having this given knowledge, relocating the idle vehicles to areas with high trip demand potential could significantly improve the serving ratio. For this relocation, Wallar et al. \cite{wallar2018vehicle} first found fixed relocation centers for the serving area based on maximum waiting time. They treat the road graph $G=(V, E)$ as a static graph, shown as Fig. \ref{fig:dynamicgraphs}(a), and formulated the relocation center searching as an Integer Linear Programming problem (ILP). Due to its high computational complexity, this method requires so long running time that can not be used for online relocation in a real time scenario. Our previous work \cite{guo2020FDAVeD} can be used online to deal with only the dynamic traffic conditions, shown as the changing attribute of edges in the road graph $G^t = (V, E, W_e^t)$ (Fig. \ref{fig:dynamicgraphs}(b)). Some traditional graph partitioning algorithms \cite{von2007tutorial, Lin2010power} could also partition the graph efficiently, but only consider the weights on edges.

In a vehicle dispatching system, road graph have dynamic traffic information on edges ($W_e^t$) and trip request information on vertices ($W_v^t$), as shown in Fig. \ref{fig:dynamicgraphs}(c). The existing research does not effectively use both these information for idle vehicle relocation. Therefore it demands further research for a vehicle dispatching system that considers the dynamic information of travel demand on the nodes and traffic on the edges, and that could search relocation centers online based on the dynamic information. The online relocation in this manner will find the most suitable relocation centres and effectively improve the serving quality from different perspectives. There are two main challenges to achieve this. The main challenge in addressing this existing limitation is to simultaneously consider both traffic and demand information for identifying appropriate vehicle relocation centres dynamically, and develop an online real-time relocation mechanism.

To address the existing limitations, in this paper, we propose a dynamic future demand aware vehicle dispatching system (called \texttt{dFDA-VeD}). It is designed for urban on-demand mobility service and can efficiently make decision for both vehicle dispatching and idle vehicle relocation. Unlike present works \cite{wallar2018vehicle, guo2020FDAVeD, von2007tutorial, Lin2010power}, \texttt{dFDA-VeD} could partition road graph and search appropriate relocation centers online based on both traffic and demand information. Overall, we make the following main contributions:

\begin{itemize}
    \item[--] We propose a vehicle dispatching system, called \texttt{dFDA-VeD}, which could dynamically search idle vehicle relocation centers online, considering the attributes on both vertices and edges of a road graph.
	\item[--] We develop a dynamic road graph partitioning based optimization objective function that considers both real-time traffic and travel demand information in order to effectively relocate the vehicles. 
	\item[--] We propose an algorithm to achieve a local optimum in order to solve the optimization objective function, and also theoretically prove its convergence. The local optimum is able approach the global optimum by parallel processing of multiple instances.
	\item[--] We perform extensive experiments on a real dataset and compare our results with the existing systems. Results show that our dynamic idle vehicle relocation based dispatching system \texttt{dFDA-VeD} outperforms the state-of-the-art vehicle dispatching systems in terms of passenger serving ratio with an increase in a very small operating cost.
\end{itemize}

The rest of the paper is organized as follows. Section  \ref{sec:relatedworks} presents the related works in vehicle dispatching system and idle vehicle relocation task. Then we give the problem formulation in Seciton \ref{sec:formulation}. The \texttt{dFDA-VeD} system is introduced in Section \ref{sec:system}.  Section \ref{sec:experiment} explore the performance of the \texttt{dFDA-VeD} system against real-world taxi data, and compared it with several baselines.  Finally we conclude the paper in Section \ref{sec:conclusion}, along with a brief discussion on the future research directions.

\vspace{-3pt}
\section{Related Work}
\label{sec:relatedworks}

Vehicle dispatching problem has been studied for decades \cite{dantzig1959truck, toth2002vehicle, pisinger2007general, kool2018attention, wang2019ridesourcing}. The objective of a dispatching system is to provide better service for the passengers, specifically higher serving ratio, shorter waiting time, lower cost and so on.  Dandl et al. \cite{dandl2018comparing} converted the dispatching problem to two bipartite matching problem: vehicle-to-user and relocation assignments. Liu et al. \cite{liu2019globally} periodically get the optimal result in offline using predicted demands and use the offline results to guide the online dispatching. The two bipartite matching problem and future demands are both considered in our previous paper \cite{guo2020FDAVeD}. Tang et al. \cite{tang2019deep} and Al et al. \cite{al2019deeppool} proposed reinforcement learning method to solve the vehicle dispatching problem. Kim et al. \cite{kim2020multi} considered multi-objective vehicle dispatching problem and use minimum cost maximum flow algorithm to solve. Liu et al. \cite{liu2020mobility} considers the mobility ride demand on roadside which are not sent to centralised platform. These researchers consider different objective and methods to improve the dispatching service quality. In this paper, we handle both the vehicle--request matching and idle vehicle relocation problem, and make several contributions on the latter one.

For a city, it is important to know the minimum number of vehicles that can serve the travel requests in the region. In 2018, the minimum number of vehicles to serve a city is addressed when all travel requests are known in advance \cite{vazifeh2018addressing}. They transferred the problem to find the minimum path cover for directed acyclic graph, where nodes stand for trips. This scenario could be treated as a special case that all idle vehicles are relocated to the passengers' pickup location in time. It shows the power of a perfect idle vehicle relocation strategy. However, in real urban on-demand mobility application, the travel demands are received in real time. In this section, we introduce the two kinds of methods to handle idle vehicle relocation problem: machine learning and other methods. 

\textbf{Machine learning methods} are used to design end-to-end machine solutions to relocate idle vehicles. Li et al. \cite{li2018dynamic} designed a reinforcement learning technique to reposition bikes in a bike-sharing system. In their methods, the whole serving area are partitioned to several cluster, and a spatio-temporal reinforcement learning model are trained to learn an optimal inner-cluster reposition policy for each cluster. Holler et al. \cite{Holler2019Deep} consider the problem from a system-centric perspective. They built a central fleet management agent to make decision for all drivers and trained policies using Deep Q-Networks \cite{mnih2015human} and Proximal Policy Optimization \cite{schulman2017proximal} algorithms. For the end-to-end machine learning solutions, the searching space is huge and hard to explain why the learned policy works. It needs a long time to train the optimal policy on simulation system before running online. However, using traditional optimization method, usually mathematically proof can be given from theory. Different from the machine learning methods, we use traditional optimization method to solve the idle vehicle relocation problem and give mathematical proof to guarantee the local optimal for the optimization.

\textbf{Optimization and heuristic methods} treat the idle vehicle relocation problem as an optimization problem and use traditional or heuristic methods to find the optimal solution. To relocate the idle vehicle, the relocation centers and subareas should be searched and idle vehicles should be redistributed the vehicles between subareas. Volkov et al. \cite{volkov2012markov}, Dandl et al. \cite{dandl2018comparing} partitioned serving area based on pickup/dropoff points' physical distance. Wallar et al. \cite{wallar2018vehicle} divided the serving area based on travel time and search the minimum number of relocation centers using linear programming methods. This method takes a long time to calculate the optimal results and can only be used in offline. Guo et al. \cite{guo2020FDAVeD} proposed a heuristic method to search the relocation centers based on the traffic conditions, which could find relocation centers in an efficient way. However, they do not consider effect of the travel demand distribution.

Our work is different from the existing idle vehicle relocation in two aspects:
(1) The existing works only consider the static information \cite{volkov2012markov, dandl2018comparing, wallar2018vehicle} or dynamic traffic information \cite{guo2020FDAVeD}. Unlike them, we consider both the dynamic traffic conditions and travel demands to search the relocation centers, and give an objective function to search optimal relocation centers. (2) Most previous studies decide the relocation centers offline \cite{volkov2012markov, dandl2018comparing, wallar2018vehicle, guo2020FDAVeD}, whereas we propose an online relocation center searching method that meets the online dispatching requirement, and guarantee a local optimal result for the objective function.

\vspace{-10pt}
\section{Problem Formulation}
\label{sec:formulation}
\subsection{Preliminaries}

The vehicle dispatching system contains two entities:  \textit{passengers} and  \textit{vehicles}. Passengers send the trip requests to the dispatching system in a streaming fashion.  After receiving a set of trip requests in a batch time, the system needs to match these requests with available vehicles. The dispatching problem is to find the best assignment plan for vehicles to serve the maximum number of requests. Here we formulate the problem with following definitions.

\begin{definition}[\textbf{Trip Request}] A \textit{trip request} $T_i$, defined as a tuple ($t_i^p, l_i^p, l_i^d$), is a trip requested by a passenger at time point $t_i^p$ (the earliest time when the passenger can be picked up) from location $l_i^p$, to drop off at location $l_i^d$. A set of trip requests during a particular time interval (e.g., 1 minute) is denoted by $\mathcal{T} = \{T_1, T_2, ..., T_{k}\}$.
\hfill\scalebox{0.7}{$\blacksquare$}
\label{def:trip_request}
\end{definition}

\begin{definition}[\textbf{Road Graph}] A \textit{road graph} $G = (V, E, W_v^t, W_e^t)$ is a directed graph presenting the road network topology, comprising a set of vertices $V = \{v_1, v_2,...,v_{n}\}$, which are pickup/dropoff points of trip requests, connected by the set of directed edges $E = \{e_1, e_2,...,e_{m}\}$, which are paths in the actual road network. $W_v^t$ is the set of vertex attributes, which are pickup-dropoff gaps (defined later in Definition \ref{def:vgap}) of corresponding vertices at time point $t$. $W_e^t$ is the set of edges weights or attributes, which are travel times on corresponding edges at time point $t$.
\hfill\scalebox{0.7}{$\blacksquare$}
\label{def:road_graph}
\end{definition}

\begin{definition}[\textbf{Pickup--dropoff Gap}] The pickup--dropoff gap $g_i^t$ is the difference between the number of pickup and dropoff demands for point $v_i$ in the time window $[t, t+t_r)$. Here $t_r$ is the time length for relocating idle vehicles  to their destinations.
\hfill\scalebox{0.7}{$\blacksquare$}
\label{def:vgap}
\end{definition}

\begin{definition}[\textbf{Served Trip}] A trip request $T_{i}$ is called as a \textit{served trip} if the passenger is actually picked up  between $t_{i}^p$ and $t_{i}^p + \Delta$, where $\Delta$ is a pre-defined serving threshold (passenger's maximum waiting time). The set of trip requests already served by the dispatching system is denoted by $\mathcal{T}_{served}$, which is a sub set of the total set of trip requests $\mathcal{T}$. Mathematically, $\mathcal{T}_{served} = \{T_i| t_i^p \leq \, ^{a}t_i^p  \, \leq t_{i}^p + \Delta\}$, where $^{a}t_i^p$ is the actual pickup time.
\hfill\scalebox{0.7}{$\blacksquare$}
\label{def:served_trip}
\end{definition}

\begin{definition}[\textbf{Served Trip Ratio}] Given the set of served trips $\mathcal{T}_{served}$ and all the trip requests $\mathcal{T}$, the \textit{serving ratio} $R$ of the centralized vehicle dispatching system is defined as the ratio of $\mathcal{T}_{served}$ to $\mathcal{T}$, i.e., $R = \frac{|\mathcal{T}_{served}|}{|\mathcal{T}|}$.
\hfill\scalebox{0.7}{$\blacksquare$}
\label{def:serving_ratio}
\end{definition}

\vspace{-6pt}
\subsection{Problem Definition}

The problem considered in this paper is to develop a vehicle dispatching system that achieves a high serving ratio by dynamically relocating idle vehicles from over-supplied areas to under-supplied areas. 
A vehicle $v_i \in \mathcal{V}$ can serve only one trip request at one time. It can start to serve a new trip request only after it has arrived at the destination of its last trip. With a given number of vehicles $n_v$ on a road graph $G$, a set of real-time trip requests $\mathcal{T}$ in a batch, a set of historical trips $\mathcal{H}$ completed, the objective of the centralized vehicle dispatching system is to serve the maximum number of real-time trip requests $\mathcal{T}$ by dynamically relocating the idle vehicles, and thus achieve a high serving ratio $R$.
\begin{eqnarray}
\label{eq:ratio}
\begin{aligned}
&\mbox{maximize }\;\; R = \frac{|\mathcal{T}_{served}|}{|\mathcal{T}|}, \\
&\mbox{subject to }\;\; \forall t_i \in \mathcal{T}_{served}, \;\; t_i^p \leq \, ^{a}t_i^p \leq t_{i}^p + \Delta
\end{aligned}
\end{eqnarray}

\vspace{-3pt}
\section{Proposed vehicle dispatching system}
\label{sec:system}

This section presents the proposed vehicle dispatching system \texttt{dFDA-VeD}. We present the overall framework of \texttt{dFDA-VeD} in Section \ref{sec:systemframework} and the detailed method for the future-demand-aware dynamic relocation of idle vehicles in Section \ref{sec:dynamicrelocation}.

\vspace{-6pt}
\subsection{Dispatching system} \label{sec:systemframework}

\begin{figure}[!ht]
\vspace{-10pt}
\centering
\includegraphics[width=\linewidth]{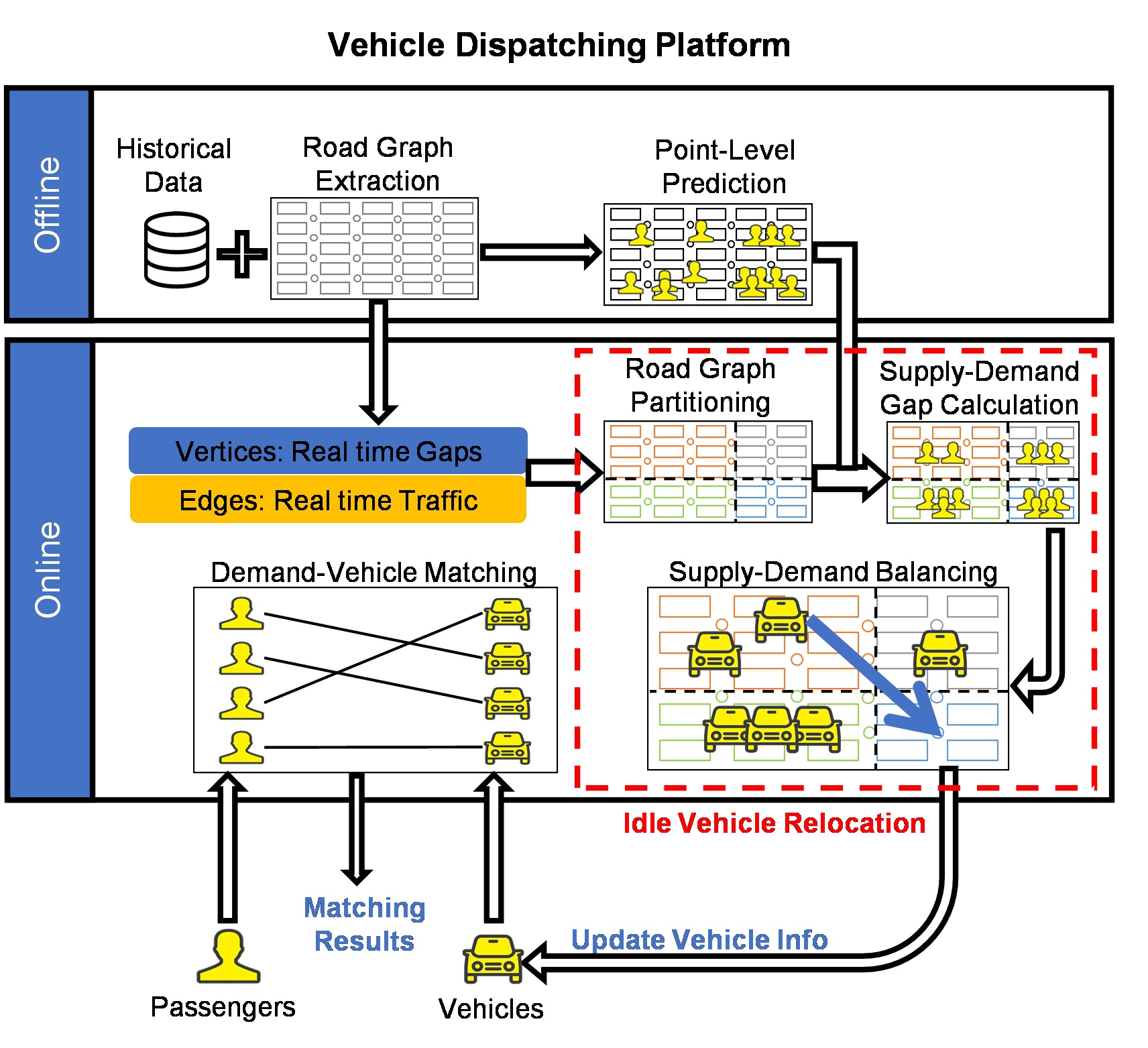}
\vspace{-10pt}
\caption{The framework of \texttt{dFDA-VeD} system}
\vspace{-10pt}
\label{fig:systems}
\end{figure}

We solve the problem of vehicle dispatching in a ride-hailing mobility-on-demand service by developing a dynamic future demand aware vehicle dispatching system (\texttt{dFDA-VeD}). Fig. \ref{fig:systems} shows the overall framework of the \texttt{dFDA-VeD} system. It starts with an offline phase of preprocessing, and then follows on to an online phase of continuously serving the realtime trip requests with available vehicles. The offline phase pre-processes the road graph data and trains a point-level travel demand prediction model based on the historical data only once in the beginning. The online phase dynamically partitions the available road graph into sub-graphs, dispatches vehicles for incoming requests, and also relocates idle vehicles based on potential future demands. The two phases of our \texttt{dFDA-VeD} system and individual modules in them are briefly explained below. As our primary focus in this paper is the dynamic relocation of idle vehicles, we skip the complete details of other tasks and modules in this paper, and present our vehicle relocation method in the next section in detail. For other tasks and modules, we follow the ideas from our previous paper \cite{guo2020FDAVeD} (refer to this paper for details).

\subsubsection{Offline Phase} There are two modules in the offline phase: extract road graph module and point-level demand prediction module. The first module extracts the road graph which is a fundamental information required in the online phase. The second module trains a point-level demand prediction model, used in the supply--demand gap calculation module (discussed later) in online phase. These two modules are briefly discussed below.

\textbf{Road graph extraction module} constructs the road graph $G=(V, E)$ for the dispatching system, where $G$ is a directed graph, $V$ is the set of vertices representing the pickup/dropoff points in the serving area, and $E$ is the set of edges representing the directed paths connecting the pickup/dropoff points in the serving area. The road graph $G$ is fundamental for all other modules.

\textbf{Point-level prediction module} uses the point-level historical average of travel demands as predictions for the prediction model. The prediction model is used later in the online phase to get predictions for the future travel demands at different pickup/dropoff points and calculate the supply-demand gap.

\subsubsection{Online Phase} The online phase is the key to dispatch vehicles and relocate idle vehicles continuously in real-time. These two tasks are performed by four modules: graph partition module, supply--demand gap calculation module, demand--vehicle matching module and supply--demand balancing module. The graph partition module dynamically partitions $G$ based on both real-time traffic and travel demand. The supply--demand gap of each subarea is calculated by the supply--demand gap calculation module. Then, the vehicle dispatching task is handled by the demand--vehicle matching module. The supply--demand balancing module is used to address the idle vehicle relocation task. These four modules are briefly discussed below.

\textbf{Road graph partition module} is used to partition the set of vertices $V$ for the whole serving area in to $k$ subareas ${A_1, A_2,..., A_k}$. The input of this module is the road graph $G=(E, V, W_e^t, W_v^t)$. Here $W_e^t$ stands for the travel time on each directed edge $e_i \in E$. $W_v^t$ stands for the pickup-dropoff gap for each vertex $v_i \in V$. For any vertex $v_i$ in the road graph $G=(E, V)$ at a specific time-point $t$, there are $pk_i$ and $dp_i$ standing for the number of pickup and dropoff demands at this vertex during the specific time interval $(t, t+\Delta)$. Then the we calculate $g_i = pk_i - dp_i$, which stands for the gap between pickup and dropoff demands in this time interval.  Here, the pickup--dropoff demand gap is the attribute of vertices $V$, so $W_v = \{g_i | i=1,...,n\}$. The subarea set $\{A_i \| i = 1,2,...,k\}$ is the partition of vertices set $V$, which means to limitations: firstly, $ \bigcup _{i=1}^k  A_i = V$; secondly, $\forall i \ne j, A_i \cap A_j = \emptyset$. The objective function for partitioning the road graph and the algorithm to optimise the objective function are presented in Section \ref{sec:dynamicrelocation}.

\textbf{Supply--demand gap calculation module} use the model trained in offline to predict the travel demand at point level, and then calculate the corresponding region level supply--demand gap for the each subarea $A_i$. The supply--demand gap is used in supply--demand balancing module to relocate the idle vehicles.

\textbf{Demand--vehicle matching module} uses Hopcroft–Karp algorithm \cite{hopcroft1973an} to find the maximum matching between the received trip requests in a short batch and available vehicles at that time. The available vehicles for a specific trip request are the vehicles that can arrive at the passenger's pickup point in $\Delta$ seconds.

\textbf{Supply--demand balancing module} relocates the idle vehicles to undersupply subareas. It starts with a search for idle vehicles. The idle vehicles are the free vehicles in over-supplied subareas. Then it follows to finding the maximum matching between the idle vehicles and the relocation centers of the under-supplied subareas. Note that the relocation centres are identified dynamically in the road graph partitioning module. The matching results are used to relocate the idle vehicles, and balance the vehicle supply in the whole serving area.


\subsection{Dynamic Idle Vehicle Relocation}
\label{sec:dynamicrelocation}

Relocating idle vehicles is very important to deal with the dynamically changing travel demand and supply of vehicles in different sub-areas of a serving area. This task would effectively re-balance the vehicles in different subareas according to the demand and supply. It requires the relocation destinations (a.k.a centres) to be firstly identified in order to make the decision. The serving area could be partitioned based on the passengers' maximum waiting time into subareas with some centres to be potential relocation destinations \cite{wallar2018vehicle, guo2020FDAVeD}. Specifically, the vehicles at relocation centers should be able to serve the trip requests in the whole serving area taking a minimum time. As the traffic conditions and the travel demands continuously change in a dynamic manner, the relocation centers also need to be updated with the changing conditions, so that they can keep serving the entire effective area in minimum time. In order to achieve this objective, we define a cost function $F(\mathcal{C}, V, W_v)$ to evaluate the performance of a set of searched relocation centers. The function is shown in Equation \ref{eq:eva_function}, where $d(c_j, v_i)$ is a distance function considered as the travel time from point $c_j$ to $v_i$, weighted by an activation function $S(\cdot)$. The distance $d(c_j, v_i)$ is the shortest travel time from vertex $c_j$ to $v_i$, and could be calculated by the attribute of edges $W_e$. The function $S(\cdot)$ transforms the  pickup-dropoff gap $g_i$ to a weight. There are multiple definitions possible for this activation function. We explore them later. The overall objective is to obtain a set of relocation centres $\mathcal{C}$ that minimise the cost function $F(\mathcal{C}, V, W_v)$, as shown in Equation \ref{eq:objective}. It is illustrated in Example \ref{ex:objective}.
\begin{eqnarray*}
\label{eq:eva_function}
F(\mathcal{C}, V, W_v) =  \operatornamewithlimits{\Sigma}_{i=1}^n \operatornamewithlimits{\min}_{c_j \in \mathcal{C}} \left\{d(c_j, v_i) S(g_i)\right\}
\end{eqnarray*}

\vspace{-8pt}
\begin{eqnarray}
\label{eq:objective}
\begin{aligned}
&\operatornamewithlimits{\mbox{minimize}}_{\mathcal{C}}\;\; F(\mathcal{C}, V, W_v), \;\;\; \mbox{subject to}\;\; 
\left\{
             \begin{array}{l}
             \mathcal{C} \subset V\\
             |\mathcal{C}| = k \\
             \end{array}
\right.
\end{aligned}
\end{eqnarray}
\vspace{-8pt}

\begin{figure}[!ht]
\vspace{-10pt}
\centering
\includegraphics[width=0.45\linewidth]{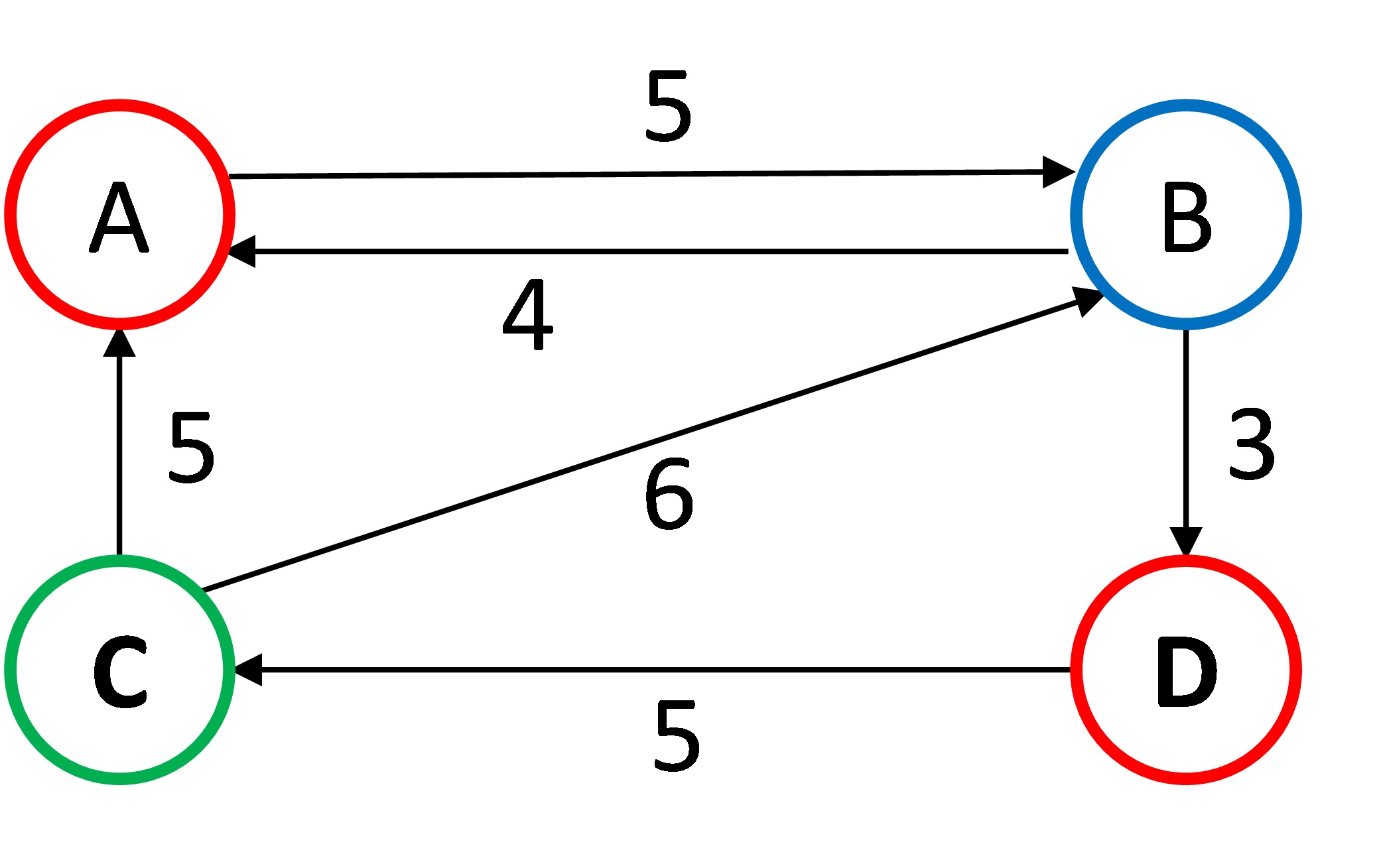}
\vspace{-10pt}
\caption{\textbf{A toy example for traffic network.} Each circle stands for a pickup/dropoff point. The red, green, and blue color stand for 1, 0, -1 of pickup--dropoff gap, respectively. The directed link stands for the way connected two points and the number show how many minutes will take through this way. For example, the directed way form A to B will take 5 minutes.}
\vspace{-10pt}
\label{fig:exampleobjective}
\end{figure}

\example{ex:objective} \textit{Consider a small traffic network of four pickup/dropoff points as shown in Fig. \ref{fig:exampleobjective}, out of which two relocation centers are to be identified. Table \ref{table:shortest_time} give the shortest travel time with and without the weighted by gap $g_i$ in destination vertex $v_i$. In Table \ref{table:objective_value}, we give the objective value with different objective. Here $F(\mathcal{C}, V) =  \operatornamewithlimits{\Sigma}_{i=1}^n \operatornamewithlimits{\min}_{c_j \in \mathcal{C}} \{d(c_j, v_i) \}$. Then with objective function $F(\mathcal{C}, V)$, the low demand vertex B and no demand vertex C would be the centers. However, with the function $F(C,V,W_v)$, the two high demand vertices A and D  minimize the function $F(C,V,W_v)$ and they would be the centers.}
\hfill\scalebox{0.7}{$\blacksquare$}

\begin{table}
    \centering
    \caption{Shortest travel time between any two vertices}
    \label{table:shortest_time}
\begin{tabular}{|l|l|l|l|l|l|l|l|l|}
\hline
  & \multicolumn{4}{l|}{d(j, i)} & \multicolumn{4}{l|}{$d(j,i)g_i$} \\ \hline
  & A    & B    & C        & D  & A(+1)  & B(-1) & C(0) & D(+1) \\ \hline
A & 0    & 5    & 13       & 8  & 0      & -5    & 0    & 8     \\ \hline
B & 4    & 0    & 8  & 3  & 4      & 0     & 0    & 3     \\ \hline
C & 5    & 6    & 0        & 9  & 5      & -6    & 0    & 9     \\ \hline
D & 10   & 11   & 5        & 0  & 10     & -11   & 0    & 0    \\ \hline
\end{tabular}
\end{table}

\begin{table}
    \centering
    \caption{Objective value with different centers}
    \label{table:objective_value}

\begin{tabular}{|l|l|l|l|l|l|l|}
\hline
         & AB & AC & AD           & BC         & BD & CD \\ \hline
$F(C,V)$   & 11 & 13 & 10           & \underline{\textbf{7}} & 9  & 11 \\ \hline
$F(C,V,W_v)$ & 3  & 2  & \underline{\textbf{-11}} & 7          & 4  & -6 \\ \hline
\end{tabular}
\end{table}

To achieve the considered objective, we develop a dynamic traffic condition and travel demand aware road graph partitioning algorithm, and use it to relocate the idle vehicles. The algorithm is based on the ideas of k-medoids \cite{Jin2010, kaufman2009finding, schubert2019faster, park2009simple}. As shown in Algorithm \ref{alg:wK-Medoids}, there are three steps to find the relocation centers and subareas. The first step initialises by randomly selecting $k$ nodes from the road graph as relocation centers (as shown in Line \ref{wK-Medoids:random_select}). The second step uses the selected relocation centers $\mathcal{C}$ as centroids to find the subareas $\mathcal{A}$ and calculate the objective value $O(\mathcal{C})$ (as shown in Line \ref{wK-Medoids:subareas}-\ref{wK-Medoids:objective}). To be simple in the algorithm, we use $O(\mathcal{C})$ stands for the objective value which means $O(\mathcal{C}) = F(C,V,W_v)$. The last step calculates new relocation centers $\mathcal{C'}$ based on the subareas $\mathcal{A}$. If $O(\mathcal{C'}) < O(\mathcal{C})$, then second and third steps are repeated until $O(\mathcal{C'}) >= O(\mathcal{C})$. Until the objective value does not decrease, we can get the relocation centers $\mathcal{C}$ and subareas $\mathcal{A}$ (as shown in Line \ref{wK-Medoids:recurrent_s}-\ref{wK-Medoids:recurrent_e}).

\begin{algorithm}[h!]
\SetAlgoVlined
\small
\KwIn{$d(\cdot, \cdot)$, $k$, $V=\{v_i | i \in (1, 2,..., n)\}$, \{$g_i$\}}
\KwOut{$\mathcal{A} = \{A_i | i \in (1, 2,..., k)\}$, $\mathcal{C} = \{c_i | i \in (1, 2,..., k)\}$}

\nonl\textbf{\uppercase\expandafter{\romannumeral 1}. Randomly select $k$ nodes as the relocation centers}\\ 

$ \mathcal{C} = \{v_i| i \in random(1, n, k)\} $ \tcp*{Initialize $\mathcal{C}$} \label{wK-Medoids:random_select}
\smallbreak 
\nonl\textbf{\uppercase\expandafter{\romannumeral 2}. Find the subareas $\mathcal{A}$  and calculate the objective function value $O(\mathcal{C})$ based on relocation centers $\mathcal{C}$}\\
$\mathcal{A} = \{A_j = \varnothing | j \in (1, 2,..., k)\}$ \tcp*{Initialize $\mathcal{A}$} \label{wK-Medoids:subareas}
\For{$v_i \in V$ }
{$j = \operatornamewithlimits{\argmin}_{i}d(c_i, j) S(g_i)$ \tcp*{Find nearest $c_j$ for $v_i$}
$A_j.add(v_i)$ \tcp*{Add $v_i$ for set $A_j$}}
$O(\mathcal{C}) =  \operatornamewithlimits{\Sigma}_{v_i \in V} \operatornamewithlimits{\min}_{c_j \in \mathcal{C}} d(c_j, v_i) S(g_i)$ \tcp*{Calculate the objective function value with $\mathcal{C}$}
\label{wK-Medoids:objective}

\nonl\textbf{\uppercase\expandafter{\romannumeral 3}. Search the new relocation center for each subarea, until the objective function does not decrease.}\\
\While{ True \label{wK-Medoids:recurrent_s}}
{
\For{$A_i \in \mathcal{A}$ \label{wK-Medoids:core_start}}
{
$ c_i' = \argmin_{v_i\in{V}} {\Sigma d(v_i, j)}, j \in A_i $ \tcp*{Search new center $c_i'$ for subarea $A_i$}
}
$O(\mathcal{C'}) =  \operatornamewithlimits{\Sigma}_{v_i \in V} \operatornamewithlimits{\min}_{c_j \in \mathcal{C'}} d(c_j, v_i) S(g_i)$ \tcp*{Calculate the new objective function value with $\mathcal{C'}$}
$flag = O(\mathcal{C'}) - O(\mathcal{C})$ \tcp*{calculate the difference between $O(\mathcal{C'})$ and $O(\mathcal{C})$}
$\mathcal{C}= \mathcal{C'}$ \tcp*{set $\mathcal{C}$ as $\mathcal{C'}$}
Repeat step II. to find the subareas $\mathcal{A}$ based on $\mathcal{C}$ \\
\label{wK-Medoids:core_end}
\If{$flag=0$}
{break \tcp*{If the objective function value does not decrease, then stop}}
}\label{wK-Medoids:recurrent_e}

\caption{Relocation Center Searching}
\label{alg:wK-Medoids}
\end{algorithm}

Our distance metrics are different as compared to the traditional metrics (such as Manhattan distance or euclidean metric). The distance between any two points $v_i$ and $v_j$ using traditional metrics is always the same in both ways, i.e., from $v_i$ to $v_j$ $d(v_i, v_j) \neq$ from $v_j$ to $v_i$ $d(v_j, v_i)$.

This is quite natural in a normal scenario of traffic conditions. Furthermore, the added weights on the distance measures make the calculation more complex. 
With all these calculations we need to ensure the convergence of Algorithm \ref{alg:wK-Medoids} to a minimum. We give the theoretical proof in Lemma \ref{lem:convergence} that Algorithm \ref{alg:wK-Medoids} always converges to a local minimum. To achieve a near optimal global minimum, we can run the algorithm several times with new random initial selections each time. As each run is independent to each other, the calculations of different runs can be completely parallelised easily.

\lemma{lem:convergence} \textit{The relocation centre searching algorithm (Algorithm \ref{alg:wK-Medoids}) always converges to a local minimum.}

\textsc{Proof}. Let $C$ and $A$ denote the current sets of center points and partitioned subareas respectively, and $C'$ and $A'$ denote the sets of searched new center points and partitioned subareas, respectively, by Algorithm \ref{alg:wK-Medoids}, Lines \ref{wK-Medoids:core_start}--\ref{wK-Medoids:core_end}. Two conditions need to be proved independently to guarantee that an algorithm will converge to a local minimum. \textit{First}, the algorithm should be monotonically decreasing the objective (or error) value, which means that $F(\mathcal{C}', V, W_v) \leq F(\mathcal{C}, V, W_v)$. \textit{Second}, there has to be a lower bound for the algorithm, which means that for any set of centre points $C$, $F(\mathcal{C}, V, W_v)>=\epsilon$. 

\textit{Proof of first condition}: $F(\mathcal{C}', V, W_v) \leq F(\mathcal{C}, V, W_v)$.
\begin{eqnarray}
    F(\mathcal{C}, V, W_v) &=& \operatornamewithlimits{\Sigma}_{i=1}^n \operatornamewithlimits{\min}_{c_j \in \mathcal{C}} \{d(c_j, v_i)S(g_i)\} \label{eqn:line1}\\
    &=&  \operatornamewithlimits{\Sigma}_{j=1}^k \operatornamewithlimits{\Sigma}_{v_i \in A_j}  \{d(c_j, v_i) S(g_i)\} \label{eqn:line2}\\
    &\geq&  \operatornamewithlimits{\Sigma}_{j=1}^k \operatornamewithlimits{\Sigma}_{v_i \in A_j}  \{d(c_j', v_i) S(g_i)\} \label{eqn:line3}\\
    &\geq&  \operatornamewithlimits{\Sigma}_{i=1}^n \operatornamewithlimits{\min}_{c_j' \in \mathcal{C}'}\{d(c_j', v_i) S(g_i)\}  \label{eqn:line4} \\
    &=& F(\mathcal{C}', V, W_v)
\end{eqnarray}

Here we illustrate why the inequalities in Equations \ref{eqn:line3} and \ref{eqn:line4} are valid. For Equation \ref{eqn:line3}, the subarea set $\mathcal{A}$ is fixed. For each subarea $A_i$, we search $c_i' = \argmin_{v_i\in{V}} {\Sigma d(v_i, j)}, j \in A_i$,  and this means that $c_i'$ is the point to make the function $\operatornamewithlimits{\Sigma}_{v_i \in A_i}  \{d(\cdot, v_i) S(g_i)\}$ to achieve the minimal value. So, $\operatornamewithlimits{\Sigma}_{j=1}^k \operatornamewithlimits{\Sigma}_{v_i \in A_j}  \{d(c_j, v_i) S(g_i)\} 
    \geq  \operatornamewithlimits{\Sigma}_{j=1}^k \operatornamewithlimits{\Sigma}_{v_i \in A_j}  \{d(c_j', v_i) S(g_i)\}$.
Similarly for Equation \ref{eqn:line4}, the new centre point set $\mathcal{C}'$ is fixed. For each point $v_i$, we search the nearest central point $c_j \in C$ to make $\operatornamewithlimits{\min}_{c_j' \in \mathcal{C}'}\{d(c_j', v_i) S(g_i)\}$ to reach the minimum. So, $\operatornamewithlimits{\Sigma}_{j=1}^k \operatornamewithlimits{\Sigma}_{v_i \in A_j}  \{d(c_j', v_i) S(g_i)\} 
    \geq  \operatornamewithlimits{\Sigma}_{i=1}^n \operatornamewithlimits{\min}_{c_j' \in \mathcal{C}'}\{d(c_j', v_i) S(g_i)\}$

\textit{Proof of second condition}: Now we prove that there is a lower bound of the objective function, i.e., $F(C, V, W_v)>=\epsilon$. Here $\epsilon$ stands for the optimal value for the objective function. The independent variable in the objective function is the centre point set $C$, where $C$ is a subset of the serving area set $V$, and $|C| = k$.  The total number of points in the serving area is $n$. There is a finite number of possibilities to select $k$ points of $C$ from $n$ points, which is equal to $C_k(n)=\frac{n!}{k!(n-k)!}$. Among all possibilities, there is always a minimum value for he objective function, so $F(C, V, W_v)>=\epsilon$.

The two conditions are separately proved in the above. They together prove that the relocation center searching algorithm guarantees to achieve a local minimum value for the objective function, and thus converges to a local minimum.
\hfill\scalebox{0.7}{$\blacksquare$}

\begin{figure}[!ht]
\vspace{-10pt}
\centering
  \subfigure[Ignore]{
  \label{fig:ignore}
  \includegraphics[width=0.18\linewidth]{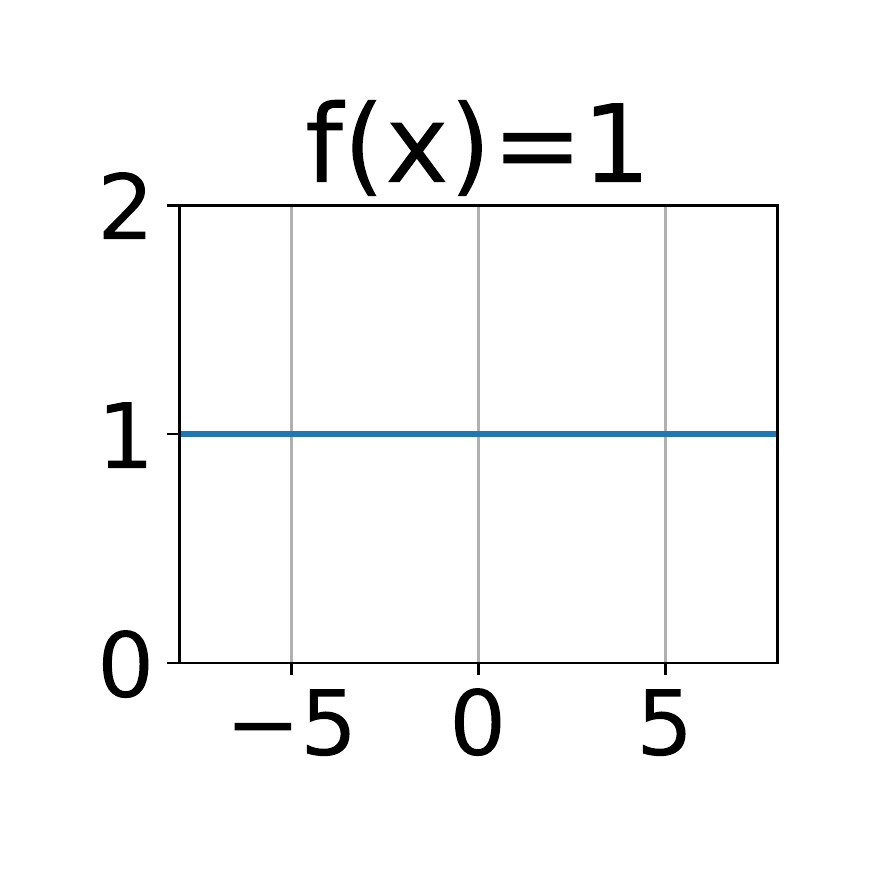}}
  \subfigure[Identity]{
  \label{fig:identity}
  \includegraphics[width=0.18\linewidth]{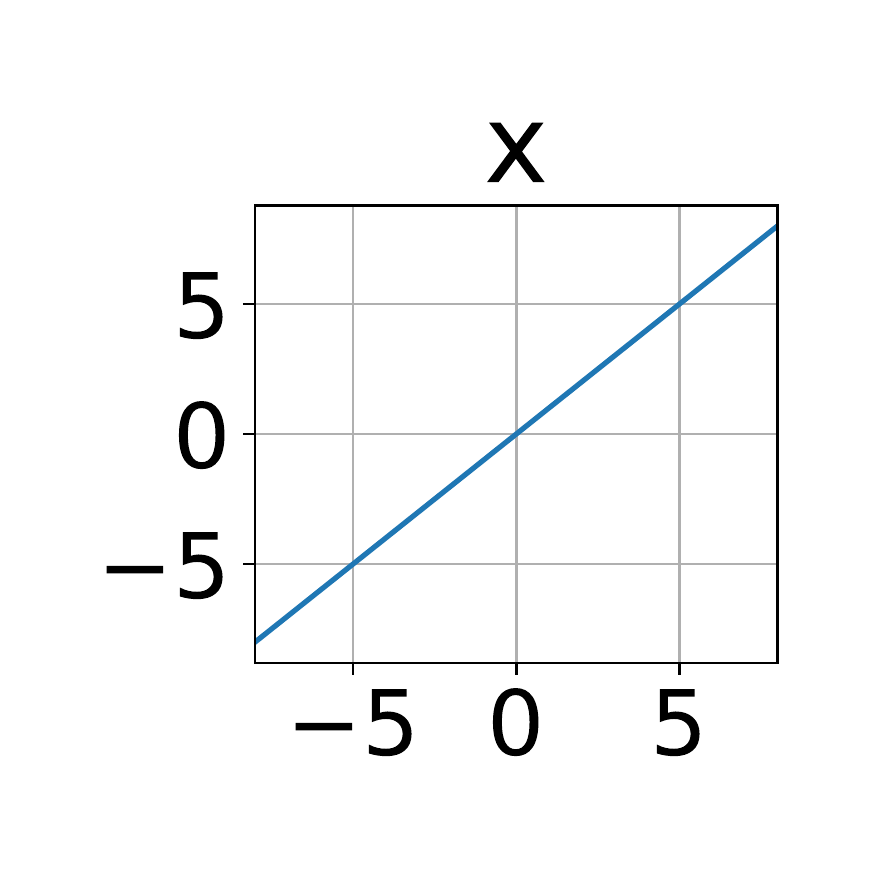}}
  \subfigure[Sigmoid]{
  \label{fig:sigmoid}
  \includegraphics[width=0.18\linewidth]{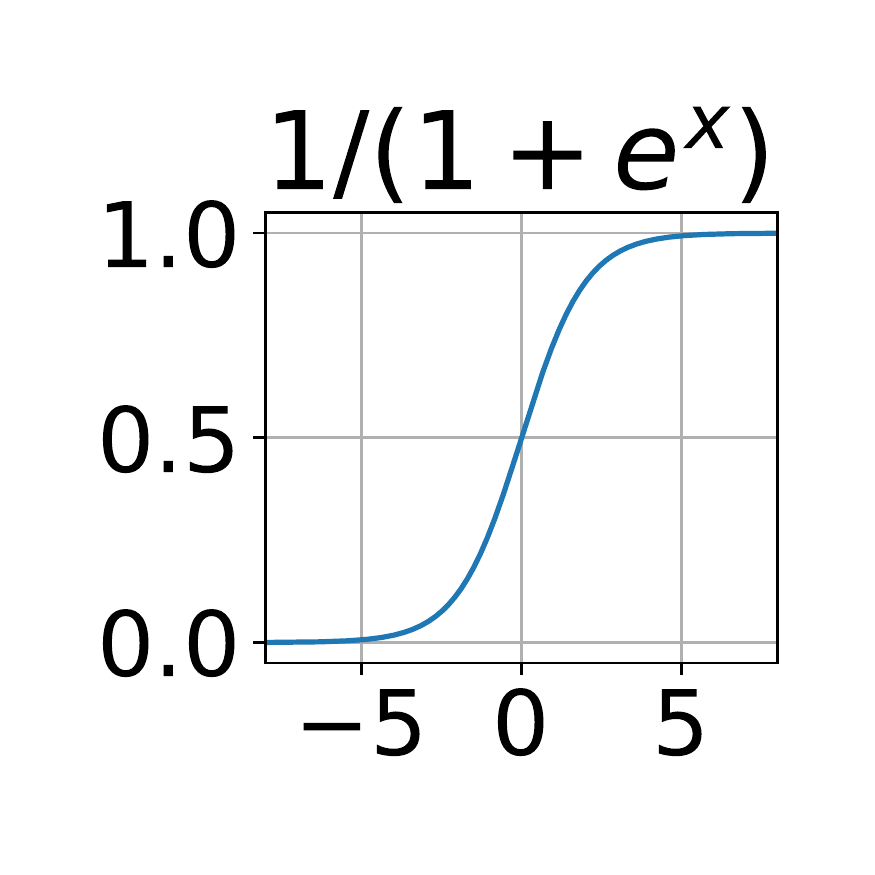}}
  \subfigure[Softplus]{
  \label{fig:softplus}
  \includegraphics[width=0.18\linewidth]{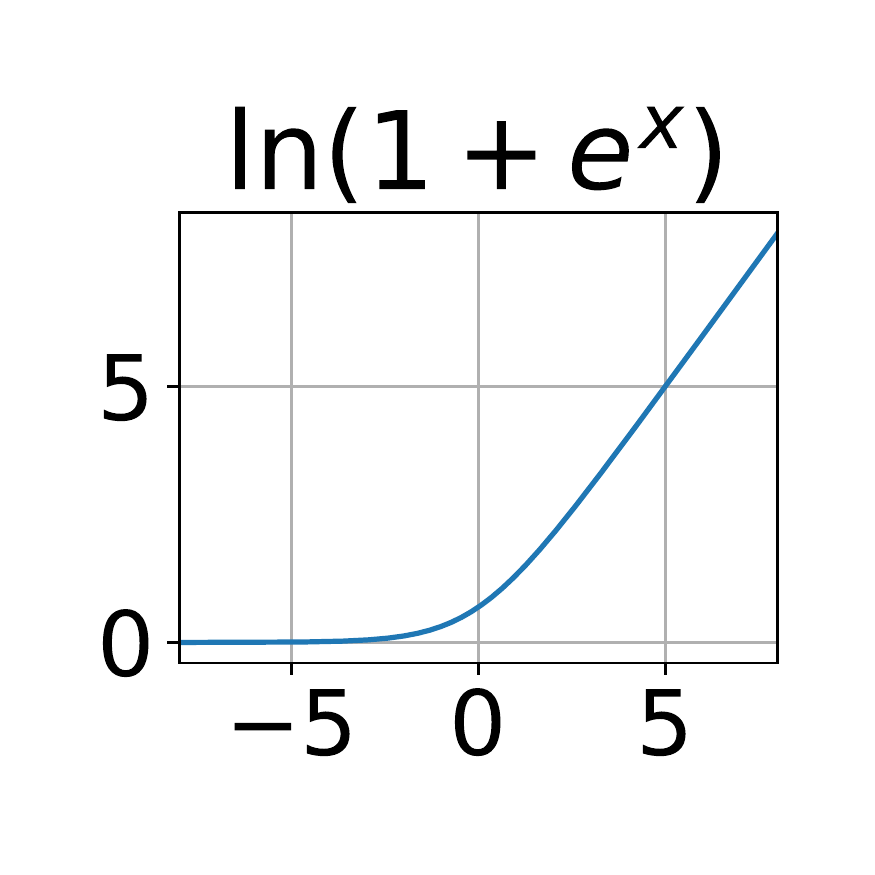}}
  \subfigure[ReLU]{
  \label{fig:relu}
  \includegraphics[width=0.18\linewidth]{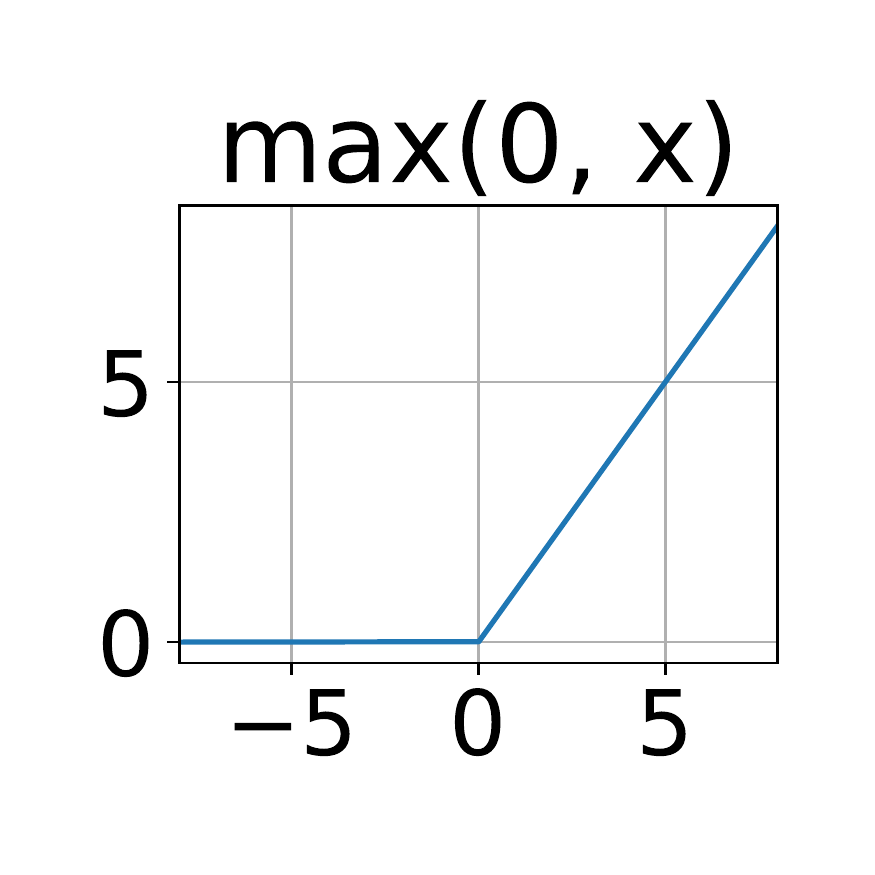}}
  \vspace{-10pt}
  \caption{Activation Functions}
  \vspace{-10pt}
  \label{fig:activations}
\end{figure}

Based on the proof, the algorithm will always converge irrespective of the exact distance metrics or the activation function. This property allow us to select several activation functions. As shown in Fig. \ref{fig:activations}, five activation functions are used in the objective function. We investigate the effects of these activation functions in Section \ref{sec:experiment}.

\vspace{-3pt}
\section{Experiments and Results}
\label{sec:experiment}

This section presents the details of our experimental evaluation of the proposed system and the obtained results. Section \ref{sec:baselines} presents the compared baselines and state-of-the-art methods. The experimental settings are described in Section \ref{sec:setting}. The evaluation metrics and experimental results are given in Section \ref{sec:metrics} and Section \ref{sec:results}, respectively.

\subsection{Compared Methods}
\label{sec:baselines}
We compare the proposed \texttt{dFDA-VeD} system having a dynamic relocation centre searching method, with four existing state-of-the-art relocation solutions: Integer Linear Programming based relocation center searching method (\textbf{ILP}) \cite{wallar2018vehicle}, spectral clustering (\textbf{SC}) \cite{von2007tutorial}, power iteration clustering (\textbf{PIC}) \cite{Lin2010power}, and our previous \textbf{FDA-VeD} system (\cite{guo2020FDAVeD}). These methods are used in our system in the road graph partitioning and relocation centers searching module for a direct comparison with the proposed method. All these methods can be categorised into three classes:
\begin{itemize}
    \item[--] The first class treats the transportation graph as a static graph (Fig. \ref{fig:dynamicgraphs}(a)). \textbf{ILP} belongs to this class.
    \item[--] The second class treats the transportation graph as a dynamic graph with time-varying edge attributes (Fig. \ref{fig:dynamicgraphs}(b)). \textbf{SC}, \textbf{PIC}, \textbf{FDA-VeD} belong to this class. 
    \item[--] The third class treats the transportation graph as a dynamic graph with time-varying node and edges attributes (Fig. \ref{fig:dynamicgraphs}(c)). The proposed \texttt{dFDA-VeD} system with novel dynamic relocation center searching algorithm belongs to this class. 
\end{itemize}

\textbf{ILP} \cite{wallar2018vehicle}: It formulated the partitioning problem as an integer linear programming problem. In this method, the maximum waiting time $\Delta$ is assumed as given, and used as a threshold. The whole serving area is partitioned into subareas to ensure that any point could be reached by its nearest relocation center less than $\Delta$. In their method, they design an offline method using a static directed graph $G=(V, E)$, a matrix $T$ where $T_{ij}$ represents the travel time between point $i$ and $j$. The road graph $G= (V, E)$ is a directed graph. The $V$ has attribute presenting the demand gap of each pickup/dorpoff location. The $E$ has attribute presenting the real-time traffic condition, specifically the travel time between the two vertices it connecting.

\textbf{SC} \cite{von2007tutorial}: We develop a spectral clustering based baseline method for partitioning the serving areas in order to identify the relocation centres. It requires a graph $G_s$ with similarity measures on the edges, for which we construct a similarity graph $G_s$ using the directed road graph $G=(V, E, W_e)$. As only edges attribute are used in this, we use $W$ to denote $W_e$ for simplicity. The weighted adjacency matrix of the directed road network $G$ is the matrix $W=\{w_{i,j}\}_{i,j=1,...,n}$, where $w_{i,j}$ stands for the travel time of an edge from vertex $v_i$ to vertex $v_j$. If $w_{i,j}=0$, it means that there is no  edge from vertex $v_i$ to vertex $v_j$. To get the adjacency matrix $A$ of the similarity graph $G_s$, firstly, we calculate the similarity of any two nodes $s_{i,j}$ using Equation \ref{eqn:sim}. In $A$, each value $A_{i,j} = s_{i, j} + s_{j, i}$. In this way, we build the similarity graph $G_s$. Then we use the standard spectral clustering method to get the partitioning.

\begin{equation} \label{eqn:sim}
s_{i,j}=
\begin{cases}
0& w_{i,j} = 0\\
1/w_{i,j}& w_{i,j} \neq 0
\end{cases}
\end{equation}

\textbf{PIC} \cite{Lin2010power}: In this method, we first construct a similarity graph $G_s$ in the same way as explained in the above \textbf{SC} method, and then partition the graph with the power iteration clustering method \cite{Lin2010power}.

\textbf{FDA-VeD} \cite{guo2020FDAVeD}: This is a method from our previous work \cite{guo2020FDAVeD}. It assumes that the maximum number of points in a subarea $n_{max}$, the total number of subareas $n_{sub} = \lfloor n_{point} / n_{max} \rfloor$ and the maximum waiting time $\Delta$ are given. Firstly it finds the points that can reach the largest number of other points in $\Delta$. Then this point is used as a central nodes and  its top $n_{max}$ nearest nodes to a new subarea. The nodes belonging to this subarea are deleted and repeat the last step to find the next subarea. If the total number of points in the serving area is N, then the partitioning result can be generated in $\mathcal{O}(n_{sub}N^2)$ running time.

In the following experiments, we investigate the these methods using the \texttt{dFDA-VeD} system and compare the results from different perspectives. 

\subsection{Experiment Settings}
\label{sec:setting}

This section presents the dataset and parameters settings of our experiments.

\subsubsection{{Dataset}} We use a trip records dataset on an actual road network, as detailed below.

\begin{itemize}
\item[--] \textbf{Trip records dataset:} The trip records dataset used in the experiments is from New York City Taxi Records in January 2011\footnote{Dataset: https://www1.nyc.gov/site/tlc/about/tlc-trip-record-data.page}. Here we select the trips whose pickup  and dropoff location are both in Manhattan Island.

\item[--] \textbf{Road Graph:} We use the road network of Manhattan Island\footnote{The road graph could be downloaded from {www.openstreetmap.org}}. Several kinds of ways (primary, primary\_link, secondary, tertiary, residential, unclassified, road and living street) are extracted to build the road graph.
\end{itemize}

\subsubsection{{Parameter Settings}}

There are several parameters need to be predefined in the modules. In the offline phase, the extract road graph module need to set the maximum length $l_{max}$ for any edge. By this way, any pickup/dropoff location could be matched to the nearest vertices in $l_{max}/2$ meters. The point-level prediction module needs two parameters: the future time interval $t_f$ and advance time interval $t_a$. These two parameters set which time length of potential demand should be predict. Here we use the historical average demands as the predicted demands. In the online phase, firstly, the total number of the vehicles should be fixed. Here we assume that all the vehicles can provide service in 24 hours.  The road graph partitioning module needs setting the number $k$ of subareas and relocation centers. The demand--vehicle matching module need the passengers' maximum waiting time $\Delta$ to find the possible vehicles to serve the request. The Supply--demand balancing module should set the idle vehicle relocation time $t_{r}$ which means that an idle vehicle should arrive the destination relocation center in $t_r$ seconds. For our experiments, we set the parameter values as follows. Maximum length of any edge $l_{max} = 200$ meters;  The advance time interval $t_a = 600$ seconds; The future time interval $t_f = 600$ seconds; The number of subareas $k$ is selected based on ILP problem's result; The number of vehicles $n_v$ is an adjusted variable in the experiments. It varies from 2,000 to 10,000 with step of 1,000; Passengers' maximum waiting time $\Delta = 300$ seconds; Idle vehicle relocation time $t_{r} = t_f$ seconds.

\subsection{Evaluation Metrics}
\label{sec:metrics}

To evaluate the performance of the proposed relocation algorithm, we deploy it into an idle vehicle relocation based vehicle dispatching system, and compare its impact to the service quality. The service quality is measured using the following metrics: passenger serving ratio, travel distance for every one kilometer, with-passenger ratio and gain--cost ratio.
\begin{itemize}

    \item[--] \textbf{Served trip ratio}, denotated as $\gamma$, is the ratio of served trips against all received travel demands.

    \item[--] \textbf{Vehicle kilometers per trip kilometer}, denotated as $\rho$.  Here, Vehicle kilometers (VKM) stands for the total travel distance of all vehicles in the dispatching fleet. Trip kilometers (TKM) stands for the total served trips' distance from pickup points to dropoff points. Then the vehicle kilometers per trip kilometer ($\rho$) could be calculated by: $\rho = \text{VKM}/\text{TKM}$.
    
    \item[--] \textbf{Trip kilometers per vehicle}, denotated as $\kappa$.  It stands for the average trip kilometers for all vehicles. Then $\kappa$ could be calculated by: $\kappa = \text{TKM}/n_v$.

    \item[--] \textbf{Trip waiting time}, denotated as $\tau$, is average trip waiting time of all served trips. The trip waiting time is the gap between the send ideal pickup time and the actual pickup time.

\end{itemize}

\subsection{{Results}}
\label{sec:results}

This section demonstrates the performance of our method in comparison to the existing state-of-the-art methods. We evaluate in terms of the quality of service in Section \ref{sec:baselinecomparsion}, in terms of operation time in Section \ref{sec:operationtime}, the impact of different activation functions in Section \ref{sec:activationfunctions}, and the road graph partitioning results in Section \ref{sec:partitionresults}.

\subsubsection{Comparison with baselines} \label{sec:baselinecomparsion}

In our experiments, we use the travel requests in Manhattan Island during 20110112-20110118 (seven days) to evaluate the impact of different algorithms to search relocation centers. For the baselines (SC, PIC and FDA-VeD)  and dFDA-VeD method, the number of centers/subareas $k$ should be predefined. Here, we use the first baseline (ILP method) to define $k$. When the passengers' maximum waiting time $\Delta=300$ seconds, the minimum number of relocation centers $k=37$. By this way, we get the number of relocation centers and regions to be partitioned. For the objective function in dFDA-VeD, we use the ReLU activation function (as shown in Fig. \ref{fig:relu}).

\begin{figure*}[!ht]
\vspace{-10pt}
\centering
  \subfigure[Served trip ratio using different algorithms for relocation center searching]{
  \label{fig:ratios_algorithms}
  \includegraphics[width=0.7\linewidth]{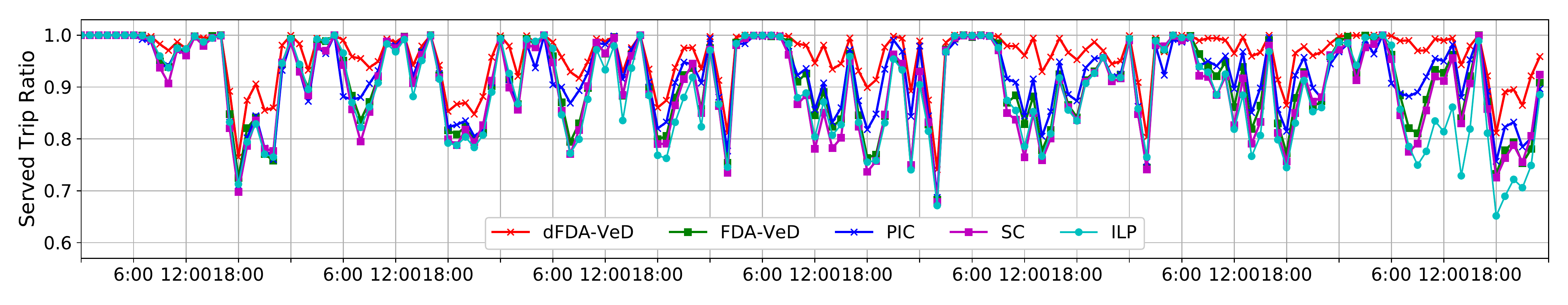}}
\vspace{-10pt}

  \subfigure[Hourly request number]{
  \label{fig:request_hourly}
  \includegraphics[width=0.7\linewidth]{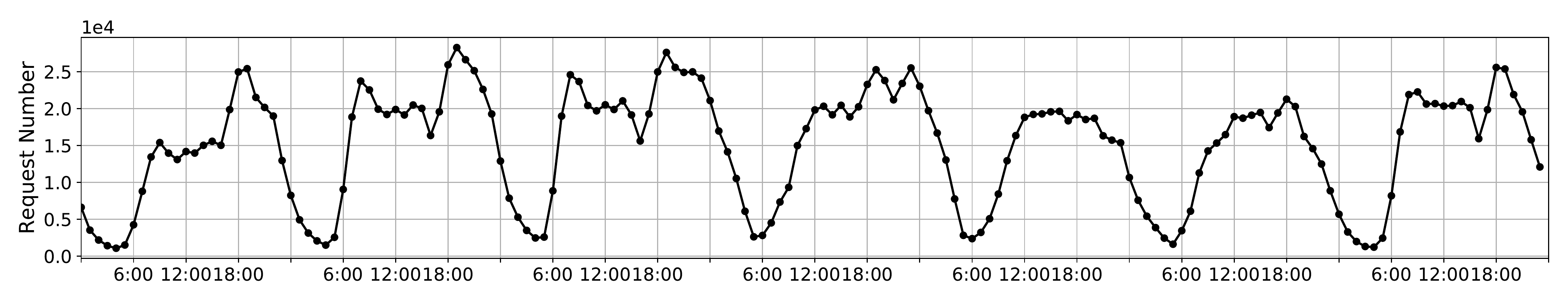}}
\vspace{-10pt}  
  \caption{Served trip ratio and request number in Manhattan Island during 20110112-20110118}
\label{fig:ratio_and_request}
\end{figure*}

\begin{figure*}[!ht]
\vspace{-10pt}
\centering
  \subfigure[Served trip ratio]{
  \label{fig:ratios_algorithms_numbers}
  \includegraphics[width=0.18\linewidth]{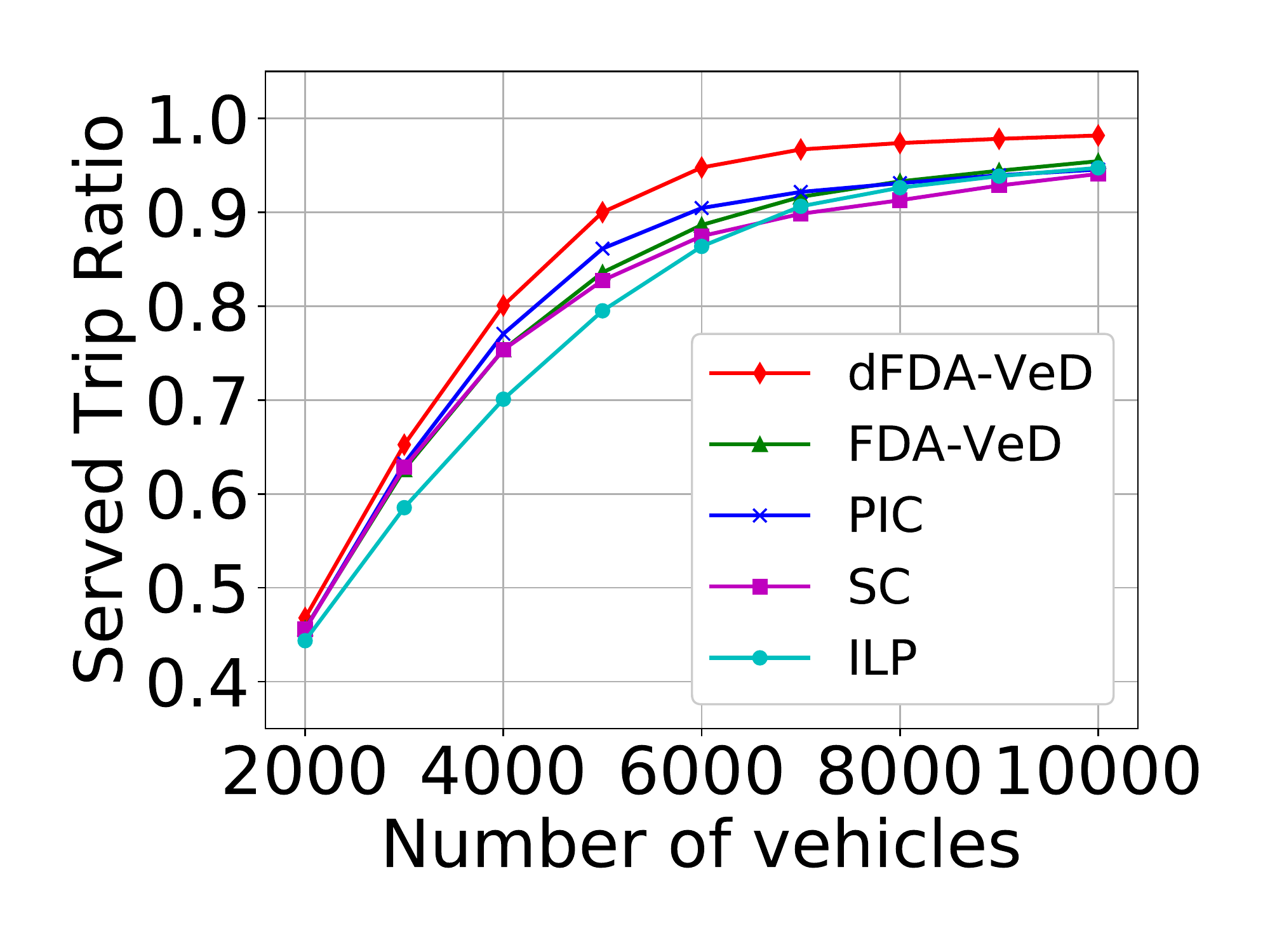}}
  \subfigure[Ratio differences]{
  \label{fig:improvents}
  \includegraphics[width=0.18\linewidth]{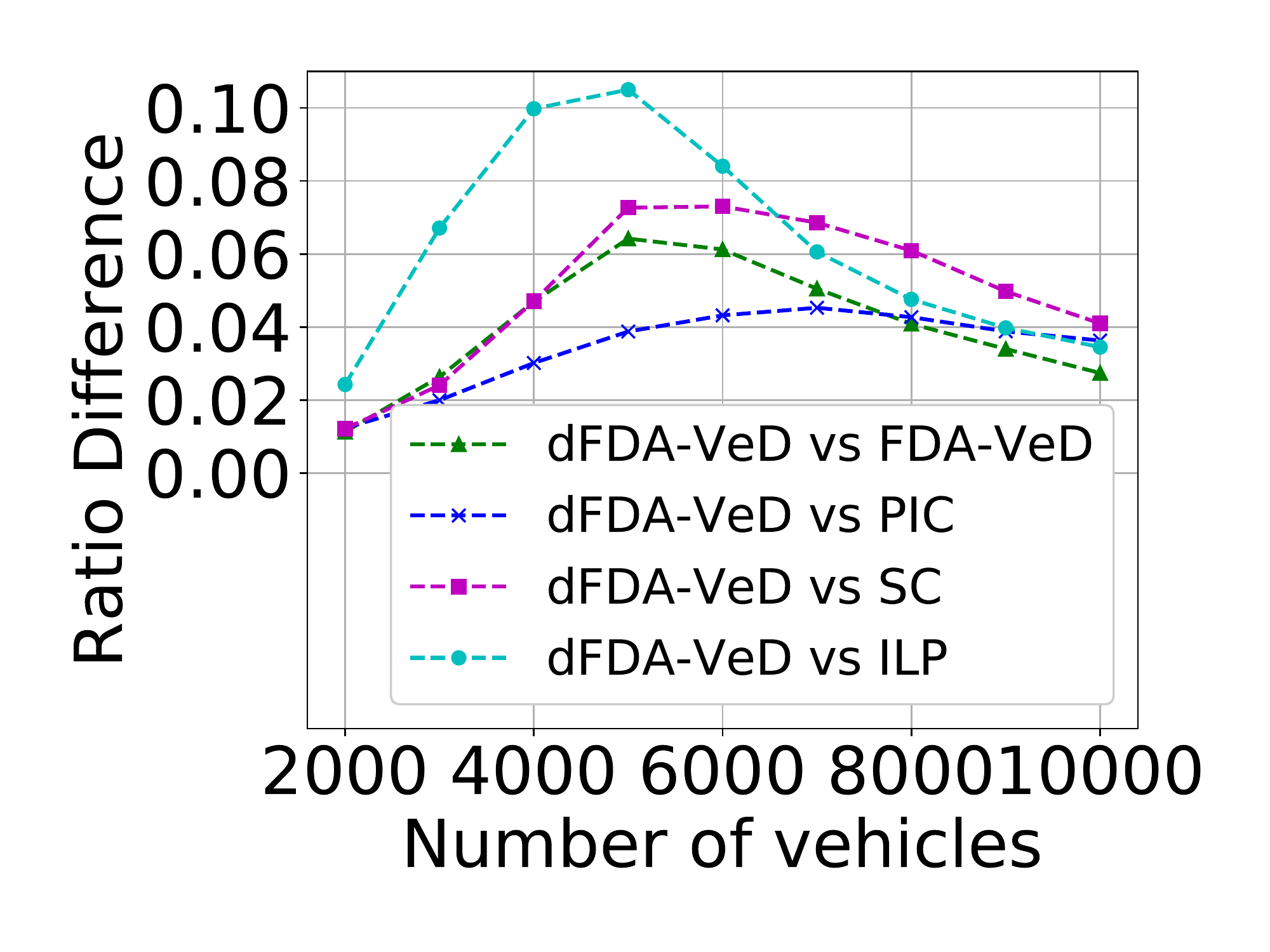}}
  \subfigure[VKM per TKM]{
  \label{fig:vkm_tkm}
  \includegraphics[width=0.18\linewidth]{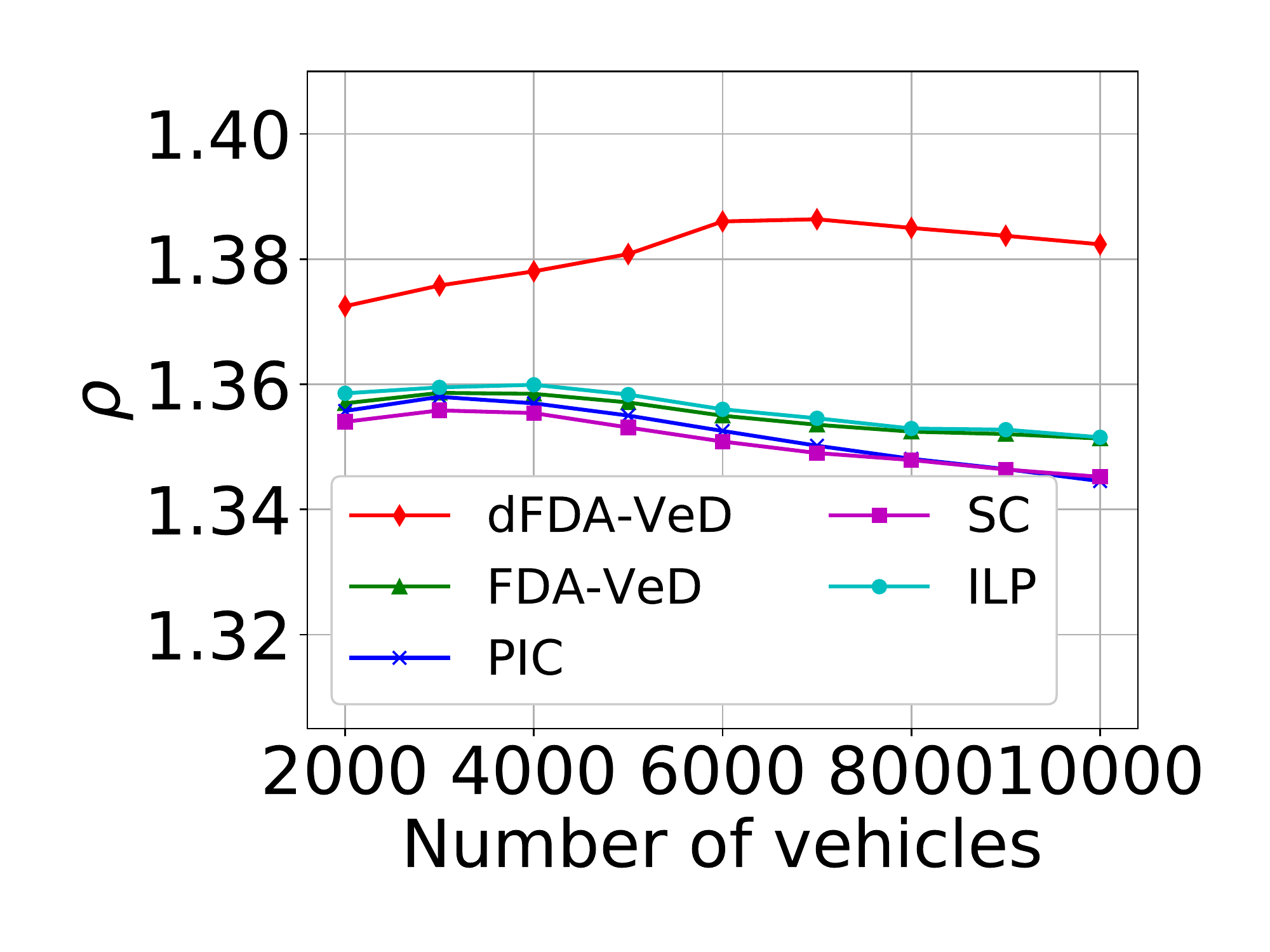}}
  \subfigure[TKM per vehicle]{
  \label{fig:trip_kilometers}
  \includegraphics[width=0.18\linewidth]{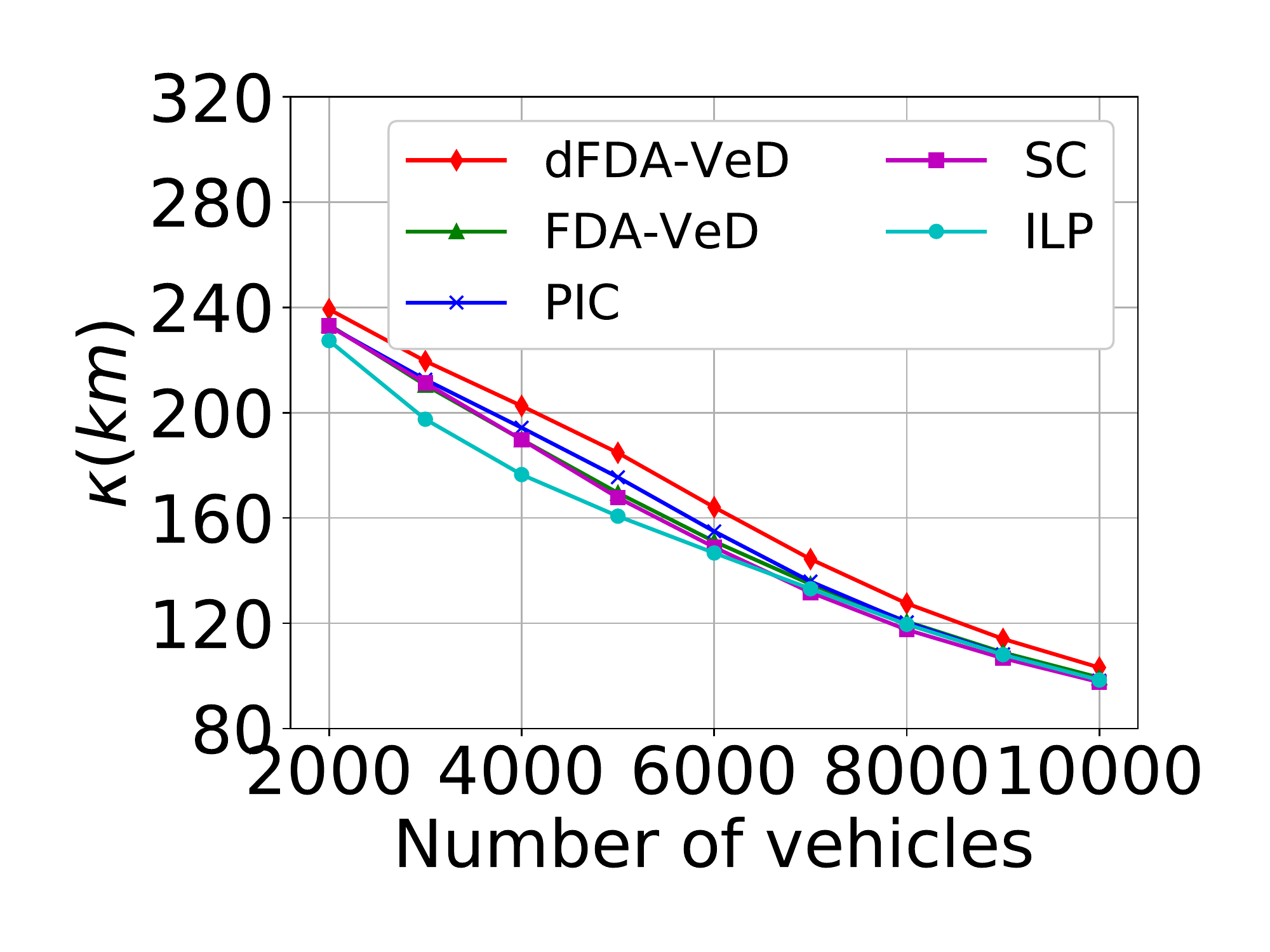}}
  \subfigure[Passenger's waiting time]{
  \label{fig:waiting_time}
  \includegraphics[width=0.18\linewidth]{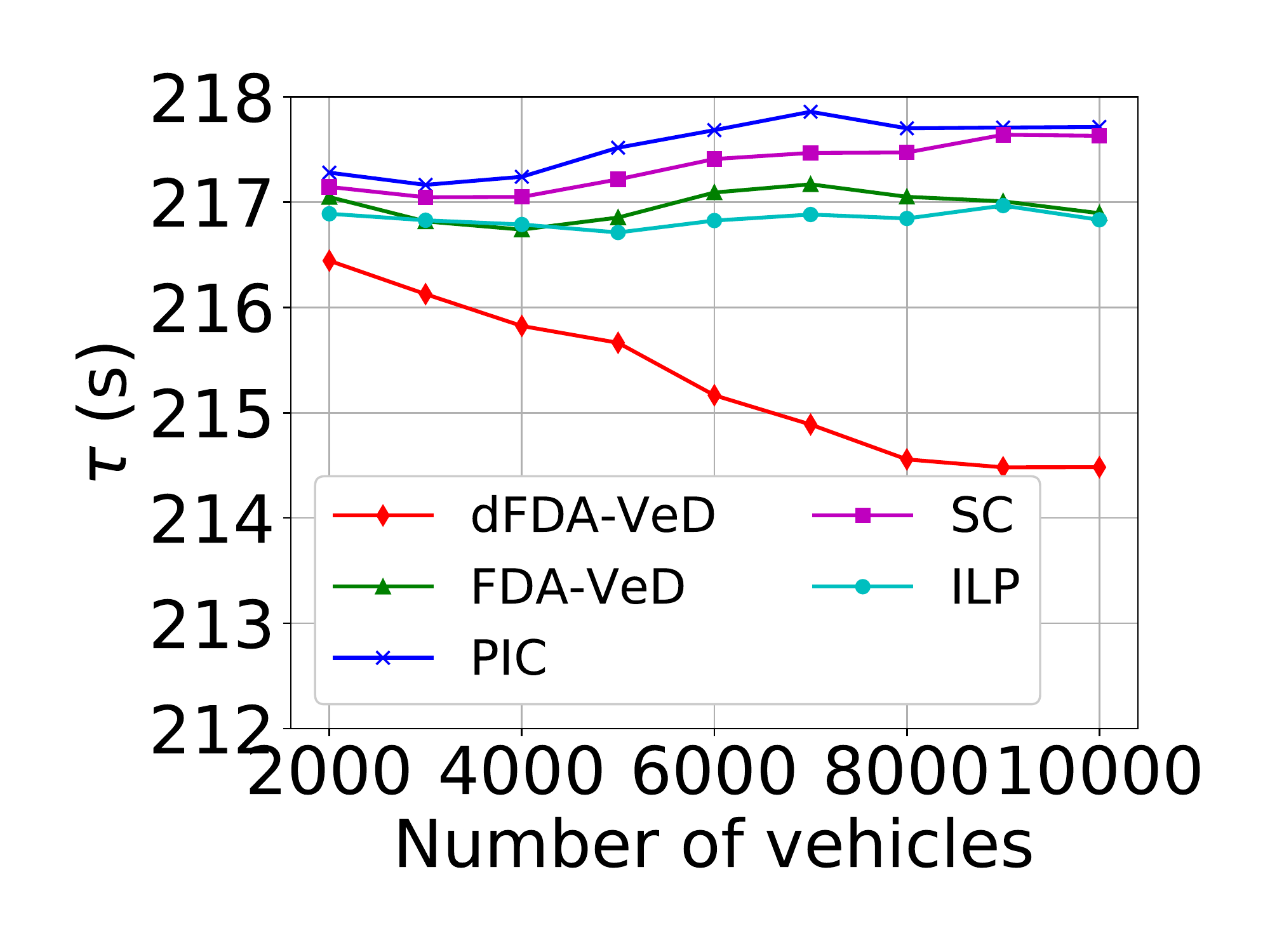}}
  \vspace{-10pt}
  
  \caption{Evaluation against different number of vehicles}
  \vspace{-10pt}
  
\label{fig:against_numbers}
\end{figure*}

\textbf{Served trip ratio ($R$):} The hourly served trip ratios $R$ using different relocation centers searching methods with 6000 vehicles are shown in Fig. \ref{fig:ratios_algorithms}.  Fig. \ref{fig:request_hourly} gives the hourly request number. Compared with the baselines, the dFDA-VeD achieved significant higher $R$, especially in peak hours (e.g. 18:00). In addition, we also investigated the served trip ratios with different number of vehicles. Fig. \ref{fig:ratios_algorithms_numbers} shows the served trip ratios in seven days when the vehicle number changes. Fig. \ref{fig:improvents} illustrates the ratio difference of the proposed dFDA-VeD against the baselines. It is clear that dFDA-VeD outperforms the baselines in all different vehicle numbers.  As shown in Fig. \ref{fig:improvents} with the vehicle number increasing, the ratio difference will raise firstly and then decline. This result is consistent with our intuition. If there is very small number of vehicles, there are very little idle vehicle to relocate. On the contrary, the baseline model could achieve a relative high $R$, and the improvement won't continue significant. The biggest gaps between dFDA-VeD and baselines appear when the number of vehicles is around 5000 or 6000.

\textbf{Vehicle kilometers per trip kilometer ($\rho$):} For this metric, the lower is better. This metric is very important for whether a dispatching could be used in real world application. If $\rho$ is very high, for example 5, it means that a vehicle drives 5 kilometers however only 1 kilometers with passengers. In other words, if $\rho$ is too high, the driver can not earn any money as the consumed oil valued more that the trip earnings. In Fig. \ref{fig:vkm_tkm}, it shows that the dFDA-VeD has a higher $\rho$ compared with the baselines. This means that there are additional cost to achieve a high served trip ratio. The increase of $\rho$ is not big. The biggest increase is less than 5\%.
    
\textbf{Trip kilometers per vehicle ($\kappa$):} For this metric, the higher is better. Fig. \ref{fig:trip_kilometers} shows that the dFDA-VeD achieves the highest $\kappa$. This means that this will bring more revenue for driver and dispatching system. 

\textbf{Trip waiting time ($\tau$):} Fig. \ref{fig:waiting_time} shows how long the passenger should wait after they send the travel requests. It shows that all these methods does not have significant difference. Although there is a downside trend of dFDA-VeD, the biggest difference of waiting time is less than 4 seconds.

\subsubsection{{Operation Time}} \label{sec:operationtime}

\begin{figure}[!ht]
\centering
  \subfigure[vehicle-request matching module (run every 1 minutes)]{
  \label{fig:time_matching}
  \includegraphics[width=0.45\linewidth]{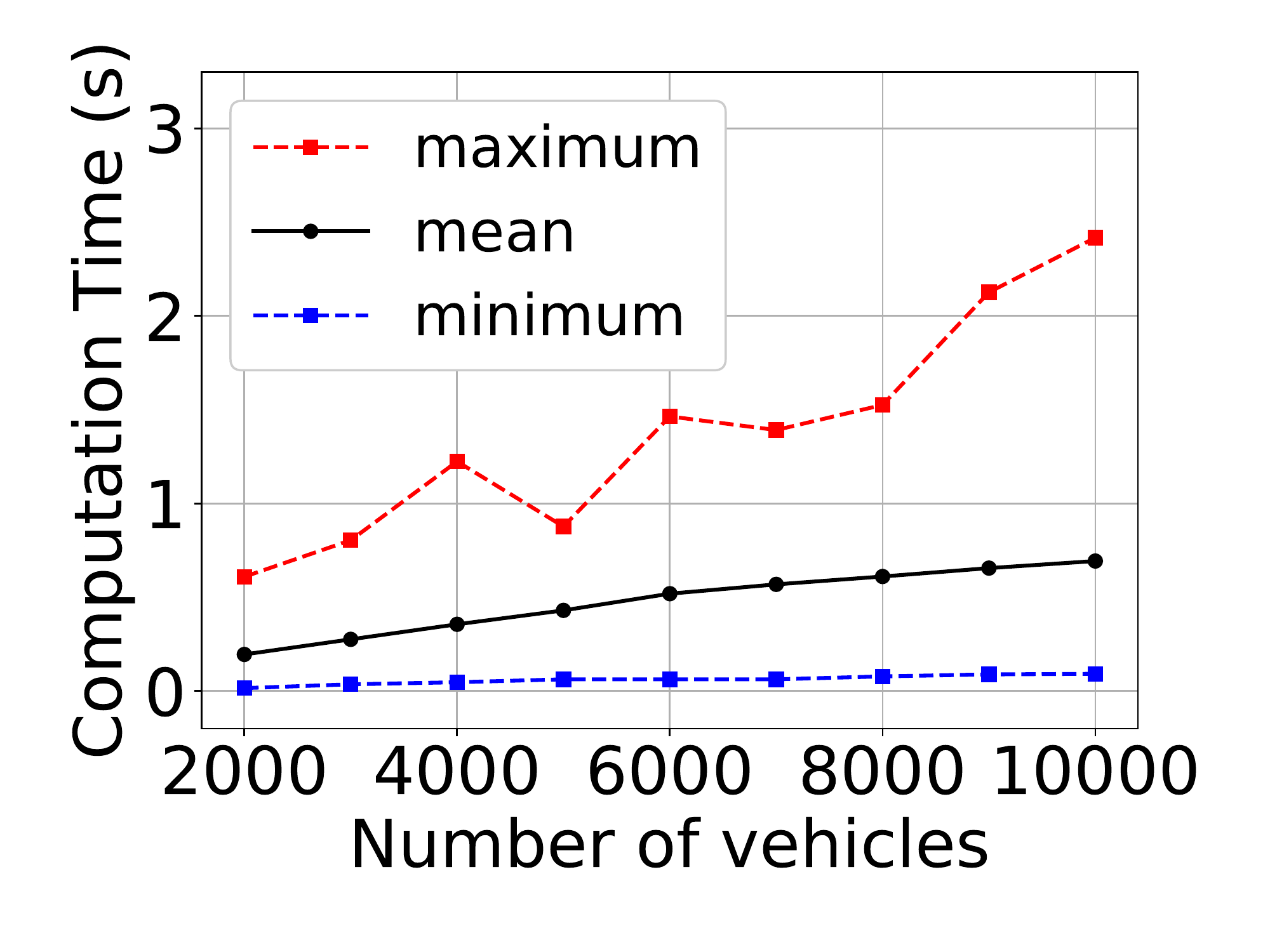}}
  \hspace{5pt}
  \subfigure[supply-demand balancing module (run every 10 minutes)]{
  \label{fig:time_balancing}
  \includegraphics[width=0.45\linewidth]{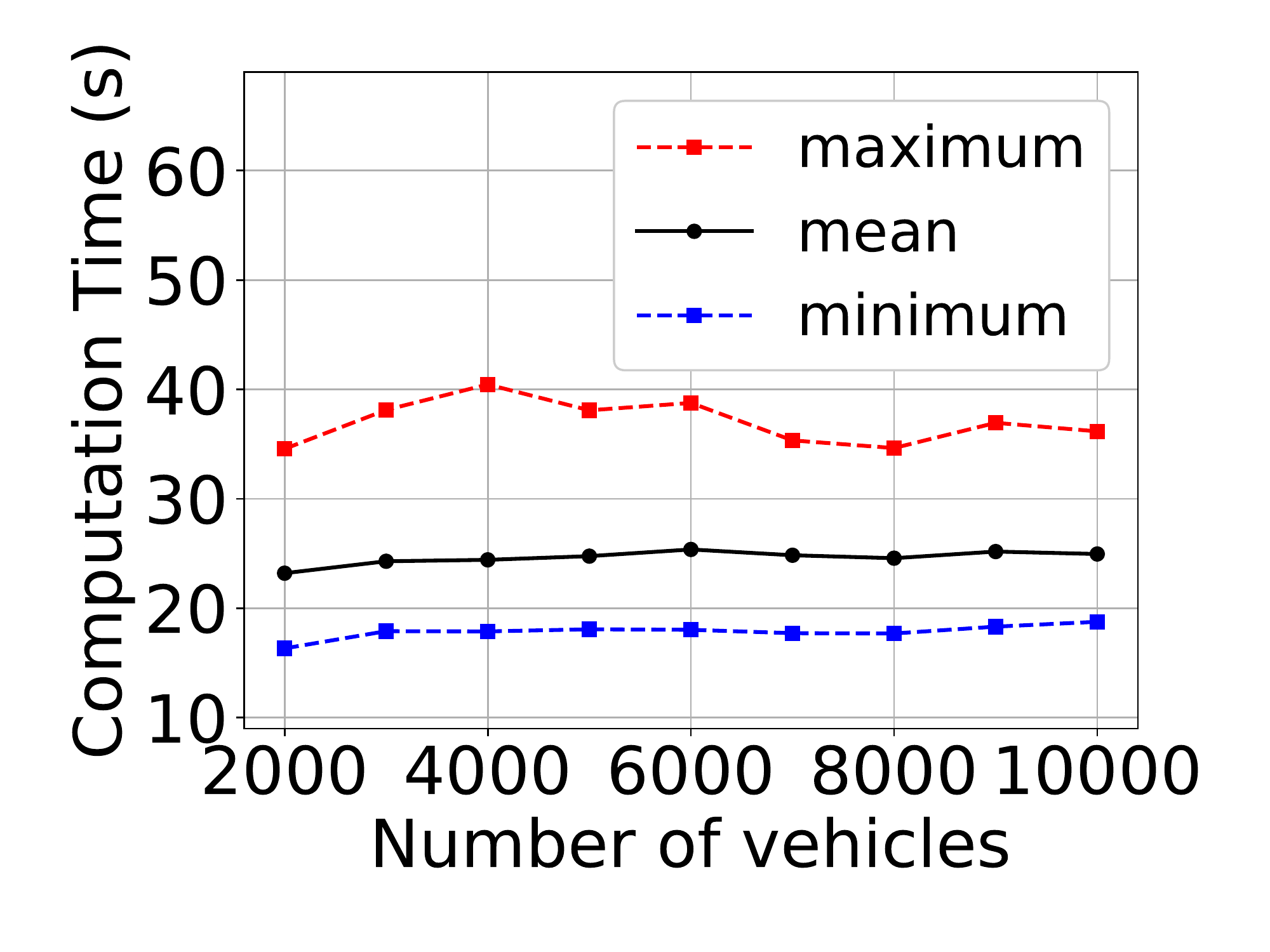}}
  \vspace{-10pt}
  \caption{Computation time for two modules}
  \vspace{-10pt}
  
\label{fig:computing_time}
\end{figure}

The  operation  time of the \texttt{dFDA-VeD} system with dFDA-VeD are shown in Fig. \ref{fig:computing_time}. Here the black lines stand for the average computation time on two task: vehicle--request matching and idle vehicle relocation. The blue areas stand for the upper bound and lower bound of the computation. For vehicle--request matching module, as shown in Fig. \ref{fig:time_matching}, it runs every 1 minutes, and the average running time is less than 1 seconds. It shows that the running time of vehicle-request matching module will go up with the vehicle number increasing.  The longest running time appears at the largest number of vehicles, as more vehicles to be matched. Even with 10,000 vehicles, the maximum computation time is less than 2.5 seconds, which can satisfy the online running requirement. For the idle vehicle relocation, as shown in Fig. \ref{fig:time_balancing}, it runs every 10 minutes, and the average running time around 25 seconds. There is no obvious trend when the vehicle number increase. For every run of idle vehicle relocation, there are three modules to run: subarea partitioning, regional level gap calculating, supply--demand balancing. The most time consuming module is subarea partitioning. At this experiment, we random select different initial central points and sequential run 8 times and select the best central points. In industry application, this time could be further decrease, as that 8 different initial points could be running in parallel. Overall, both the vehicle--request matching and idle vehicle relocation part could meet the online running requirement.

\subsubsection{Impact of activation function} \label{sec:activationfunctions}

\begin{figure*}[!ht]

\centering
\includegraphics[width=0.8\linewidth]{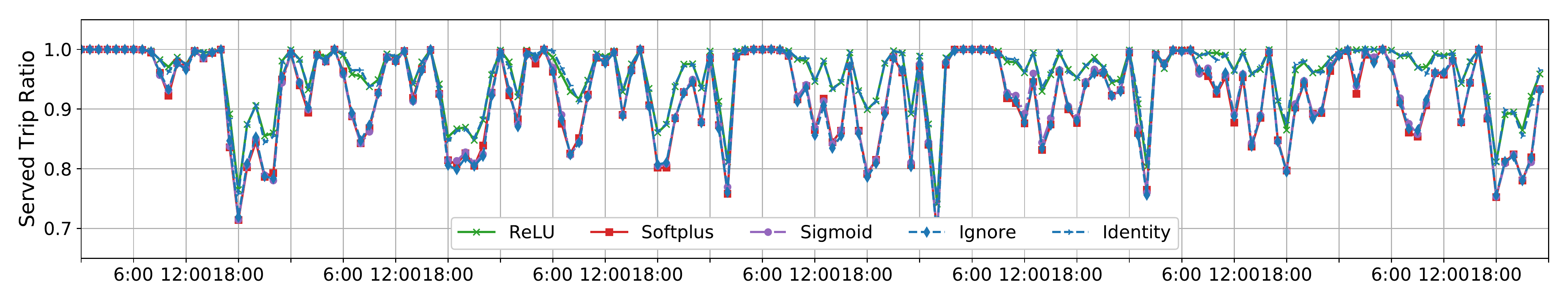}
\vspace{-10pt}
\caption{Served trip ratio using different activation functions in the optimization objective}

\label{fig:ratios_activation}
\end{figure*}

The activation function decide how the travel demand gap will affect on the distance calculation. We evaluate the served trip ratio $\gamma$ with five activation functions (\textit{Ignore, Identity, Sigmoid, Softplus} and \textit{ReLU}). Here, number of vehicles in Manhattan Island is set as 6000. Fig. \ref{fig:ratios_activation} shows that the \textit{ReLU} and \textit{Identity} activation function achieve the highest $\gamma$, these two functions achieve almost the same served trip ratio. The \textit{Ignore} function means that the travel demand gap is not considered in the optimization objective. The \textit{Sigmoid} and \textit{Softplus} functions both convert the negative gap values to positive ones, which makes the relocation centers have be more close to oversupply points than the \textit{ReLU}  and \textit{Identity} function based optimization. \textit{ReLU}  activation function transfer all negative values to zero, which means ignoring all oversupplied points when minimize the optimization objective. \textit{Identity} activation convert all under-supplied points were converted to an positive value and vice versa. The \textit{ReLU} and \textit{Identity} functions reached the same served trip ratio give us an insight: when select relocation centers, we could ignore the oversupplied points and just pay attention to the under-supplied points.

\subsubsection{Road Graph Partition Results} \label{sec:partitionresults}

\begin{figure}[!ht]

\centering
  \subfigure[The travel speed map in Manhattan Island at 18:00 (The unit of speed is m/s)]{
  \label{fig:speeds}
  \includegraphics[width=0.3\linewidth]{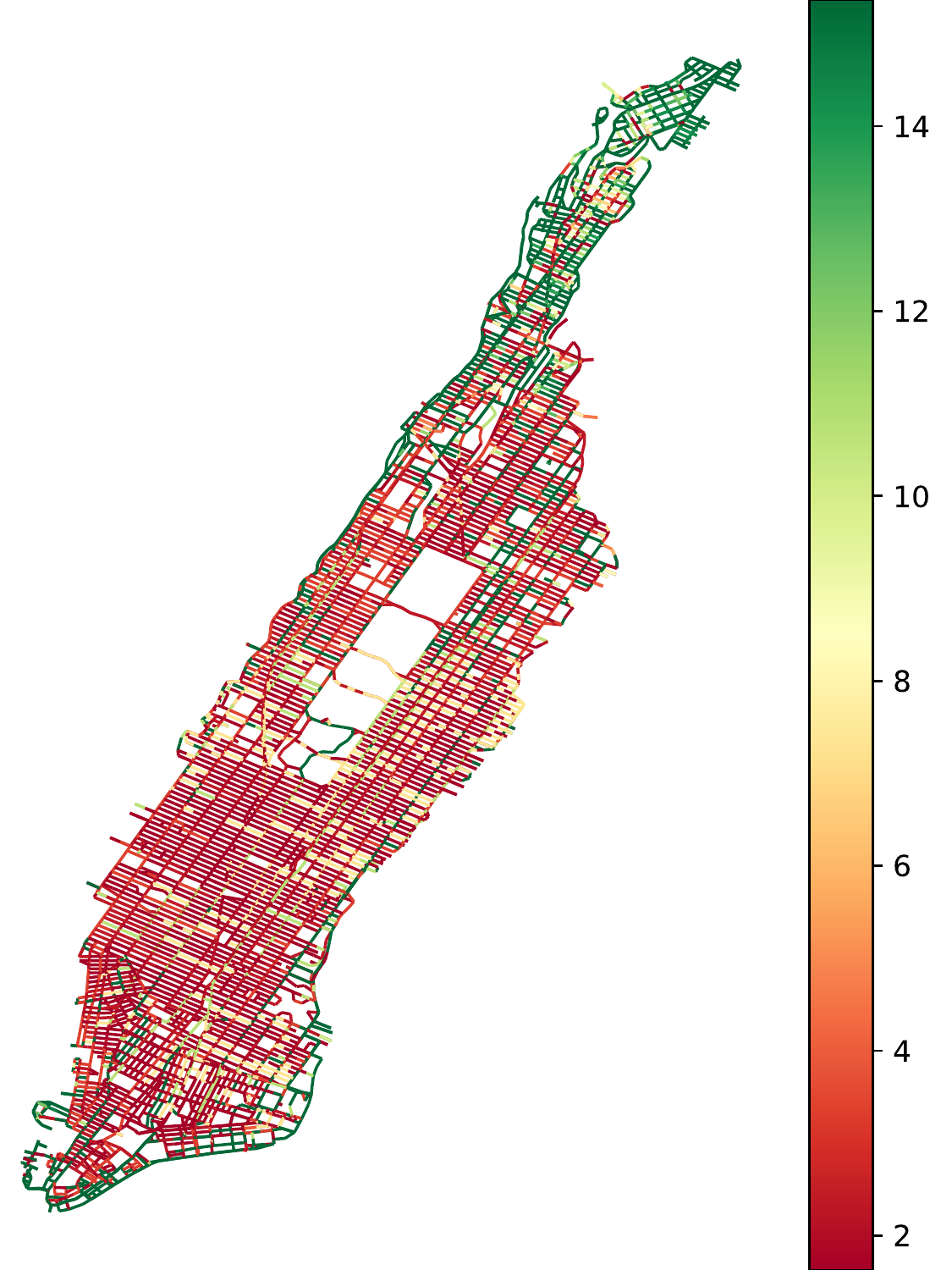}}
  \hspace{1cm}
  \subfigure[The predicted demand gap $\{g_i\}$ of each vertex during 18:00-18:10. The color of a vertex $v_i$ show the value of its demand gap $g_i$ ]{
  \label{fig:demands_gap}
  \includegraphics[width=0.3\linewidth]{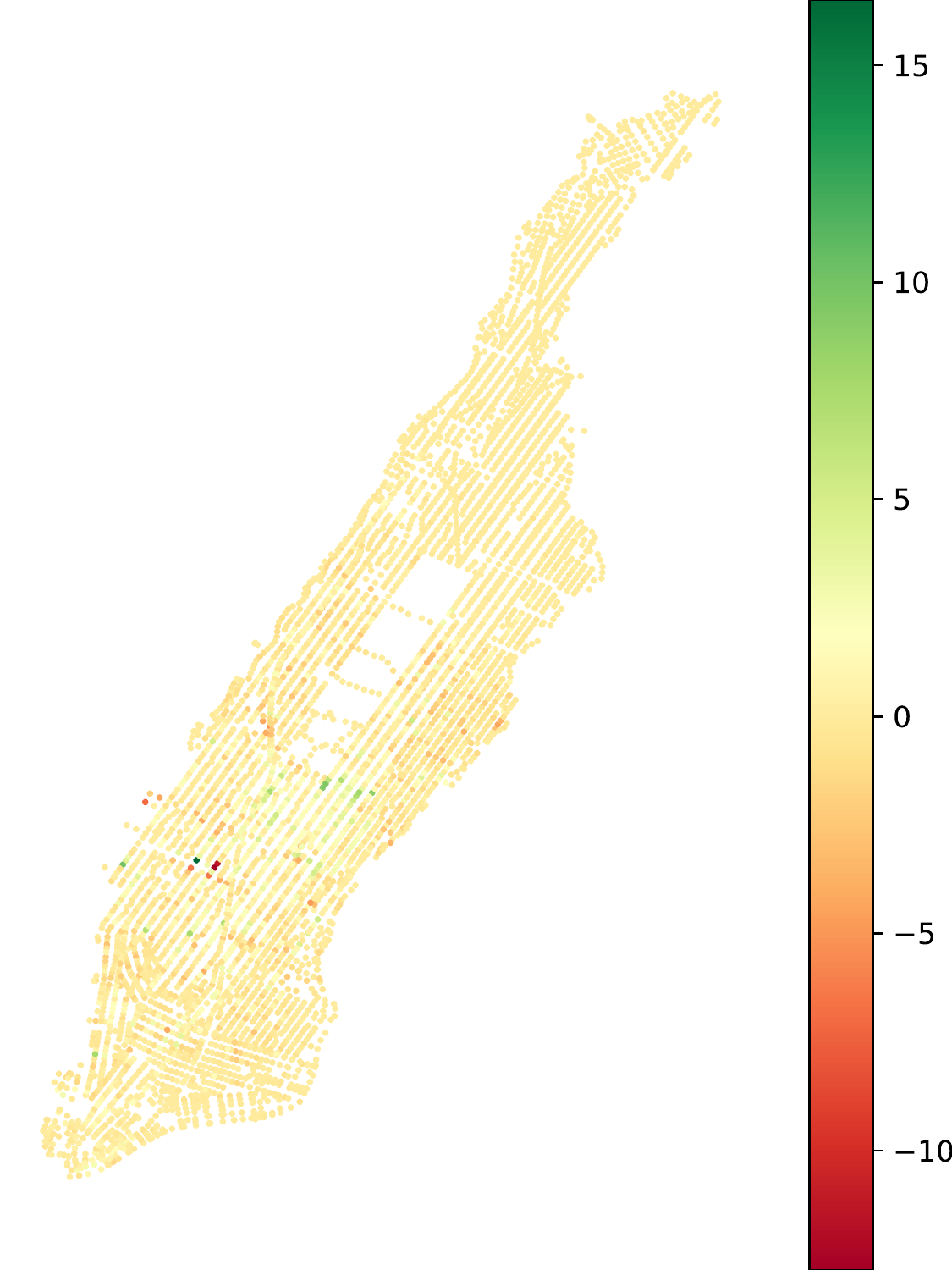}}
  \vspace{-10pt}
  \caption{The road graph of Manhattan at 18:00}
  \vspace{-10pt}
\label{fig:speeds_gaps}
\end{figure}

Here we give the road graph information (as shown in Fig. \ref{fig:speeds_gaps}) in Manhattan Island at a typical peak time 18:00. Around this time, the differences of served trip ratio between several algorithms (as shown in Fig. \ref{fig:ratios_algorithms}) and activation functions (as shown in Fig. \ref{fig:ratios_activation}) are significant. Now we would like to going deep to learn why an algorithm or activation function can beat others. Fig. \ref{fig:speeds} show travel speed on different edges of the road graph in Manhattan Island in 18:00. Fig. \ref{fig:demands_gap} show the pickup--dropoff demand gaps of each vertex during 18:00-18:10. Then with these information, we use baselines and relocation center searching algorithm in dFDA-VeD to partition the graph in sub-graphs and find the relocation centers. 

\begin{figure}[!ht]

\centering
  \subfigure[ILP]{
  \label{fig:partition_ilp}
  \includegraphics[width=0.17\linewidth]{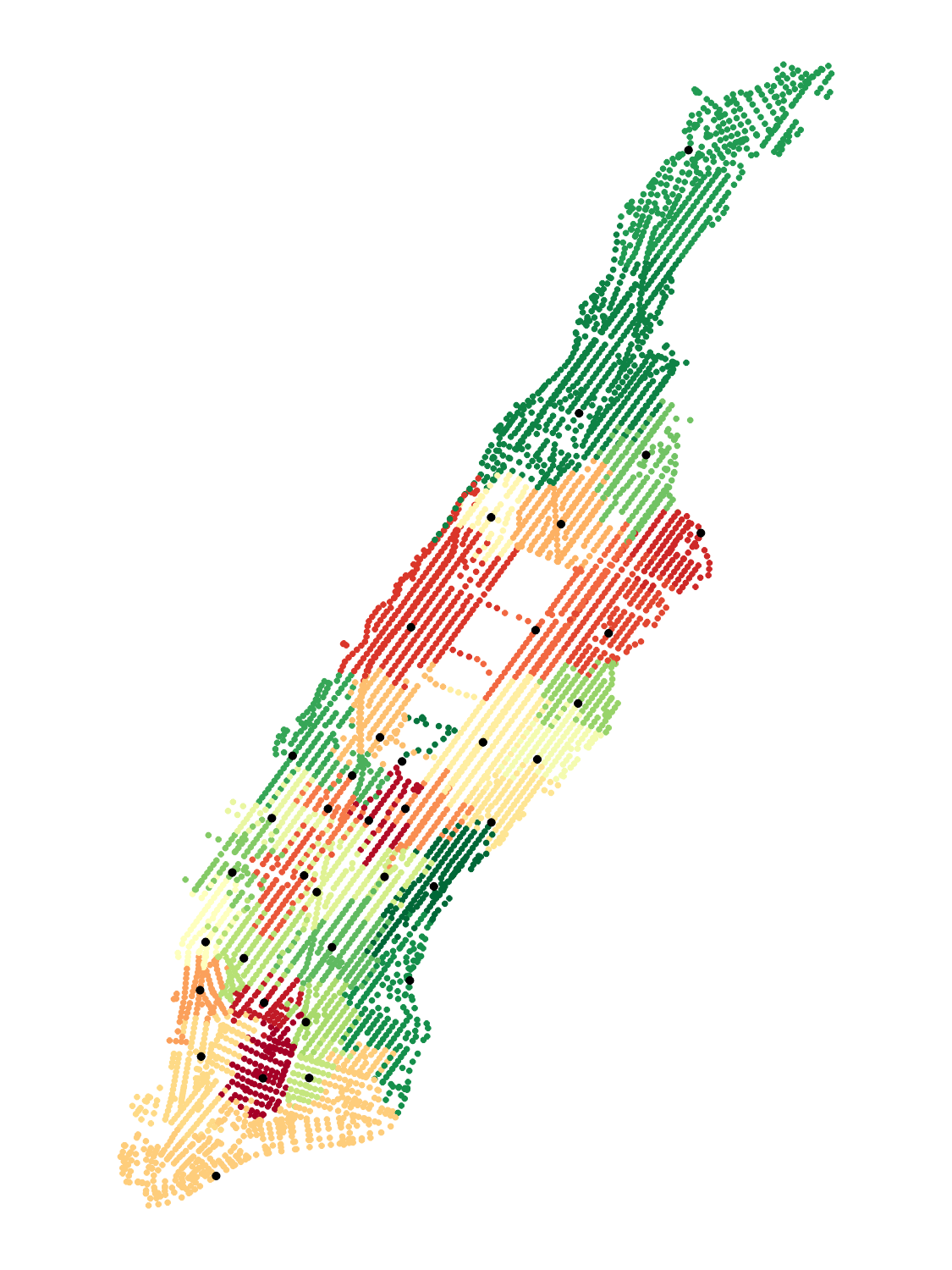}}
  \hfill
  \subfigure[SC]{
  \label{fig:partition_sc}
  \includegraphics[width=0.17\linewidth]{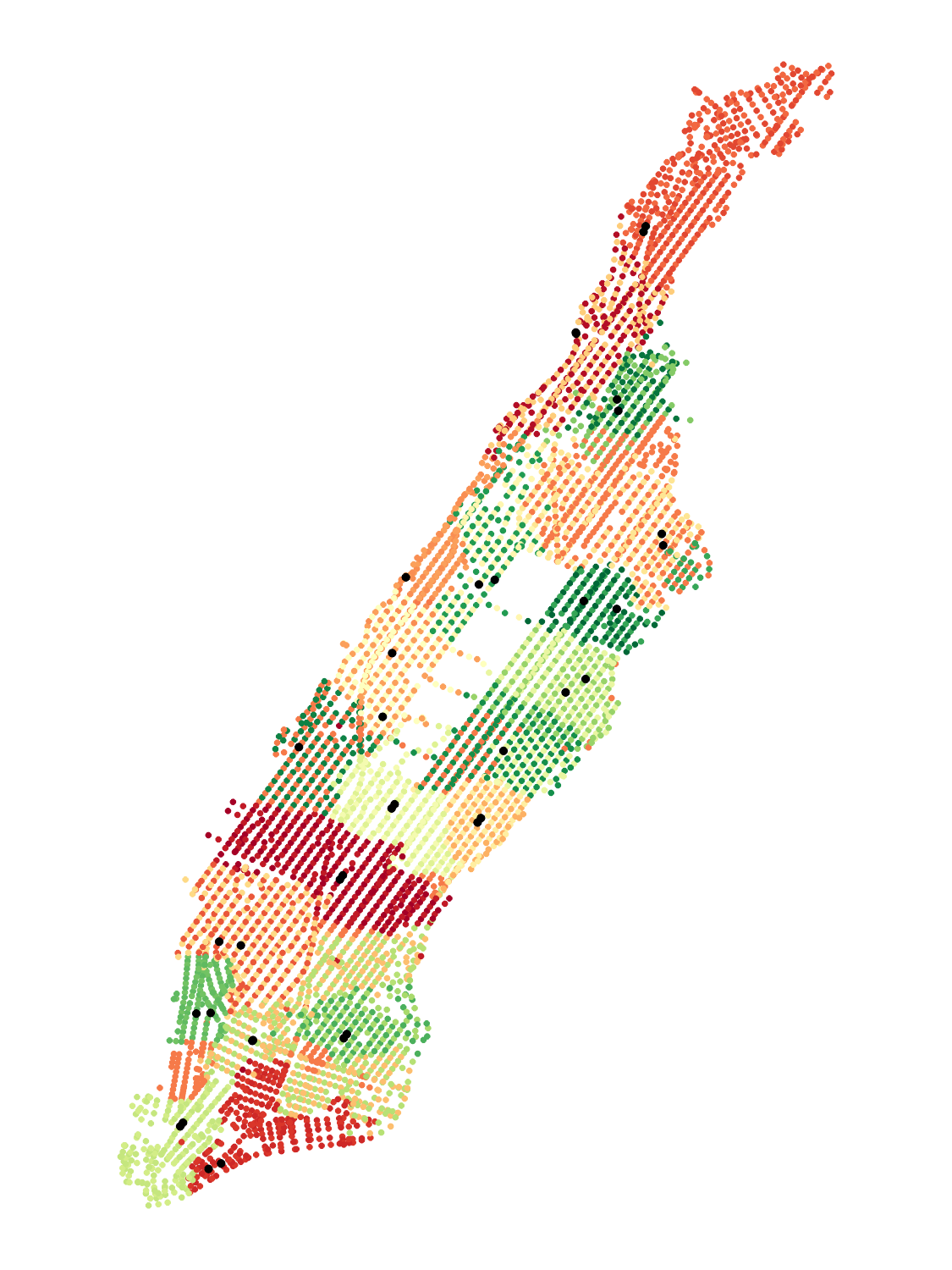}}
  \hfill
  \subfigure[PIC]{
  \label{fig:partition_pic}
  \includegraphics[width=0.17\linewidth]{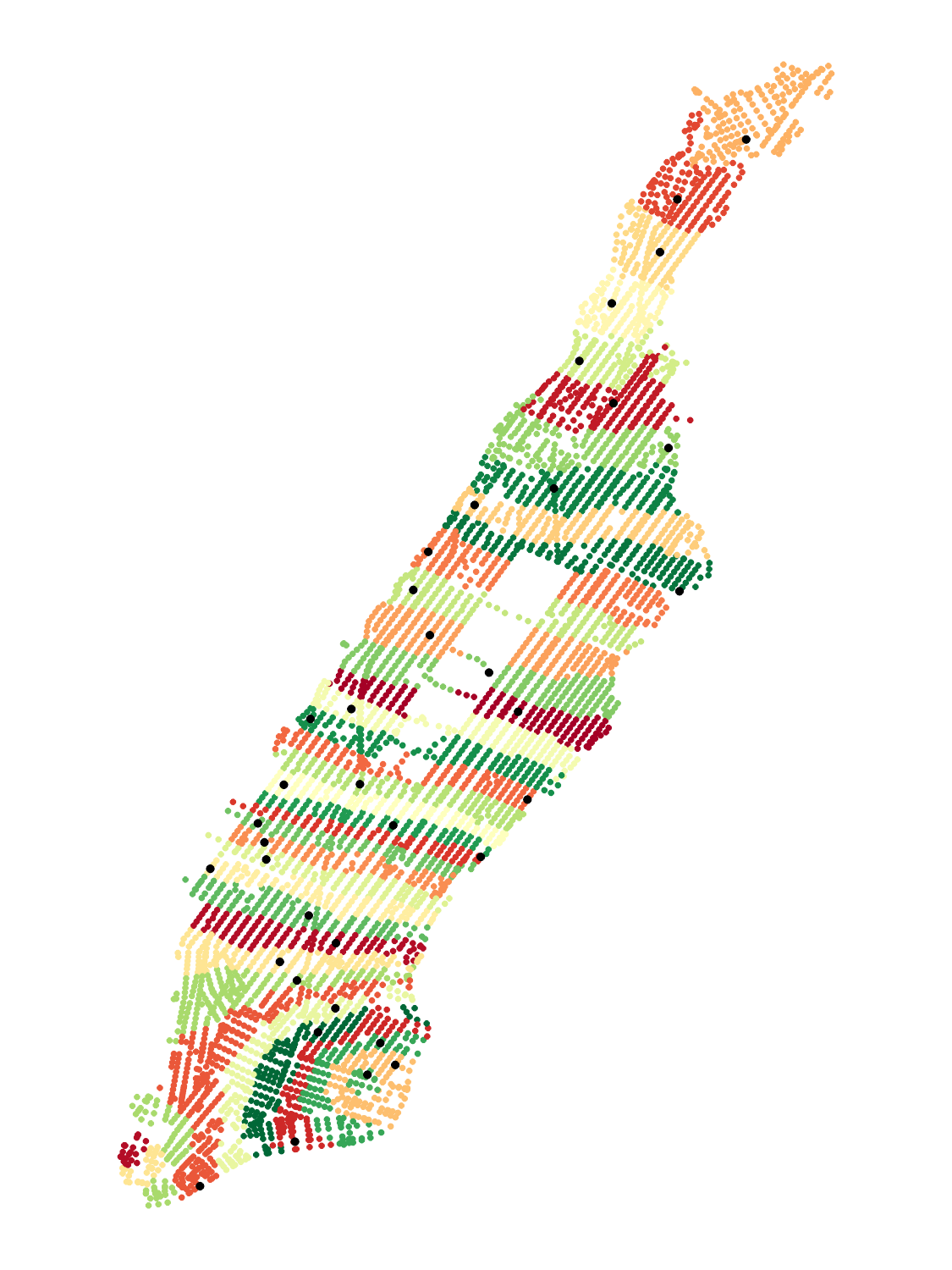}}
  \hfill
  \subfigure[FDA-VeD]{
  \label{fig:partition_fda}
  \includegraphics[width=0.17\linewidth]{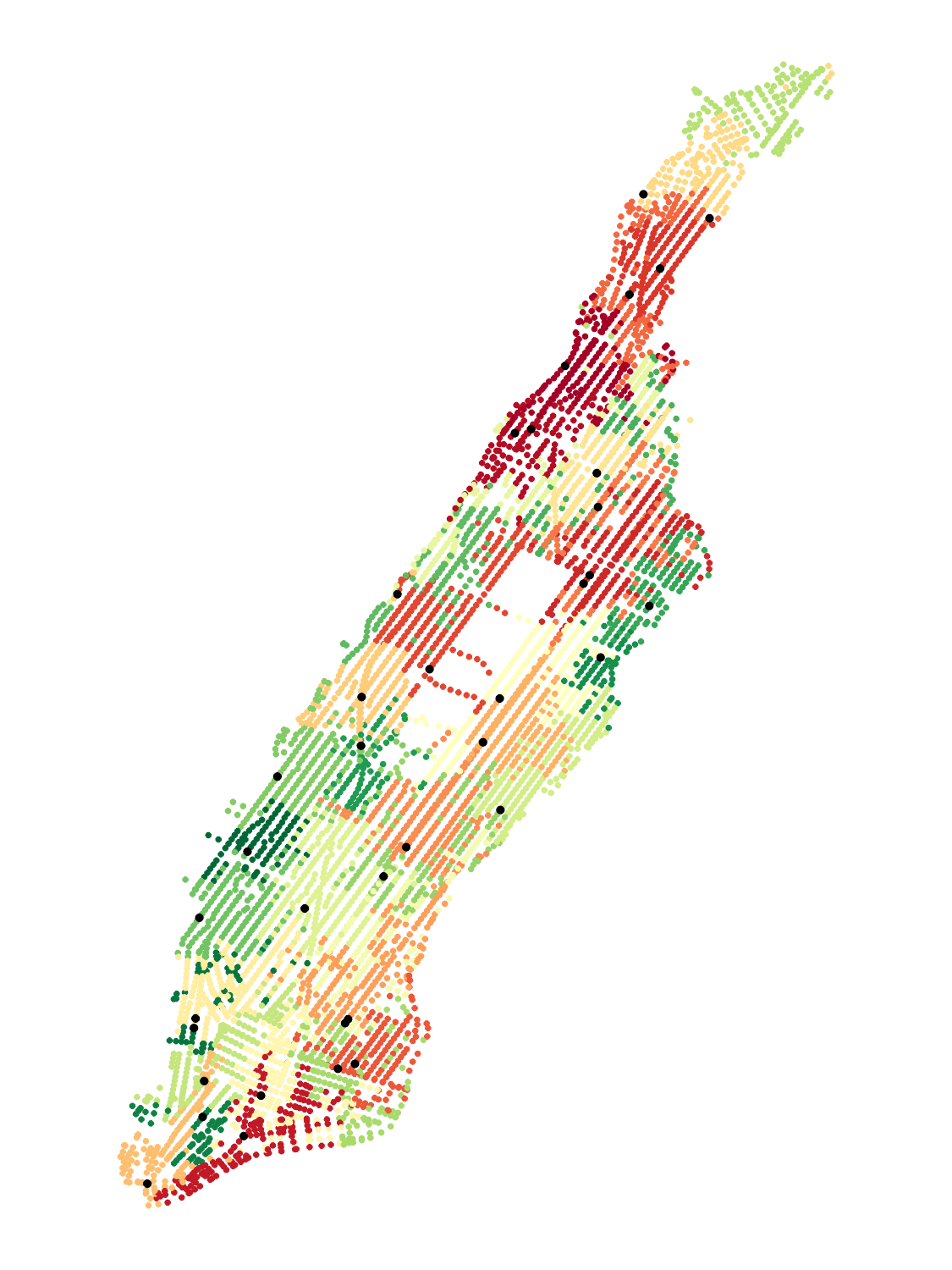}}
  \hfill
  \subfigure[dFDA-VeD]{
  \label{fig:partition_wkm}
  \includegraphics[width=0.17\linewidth]{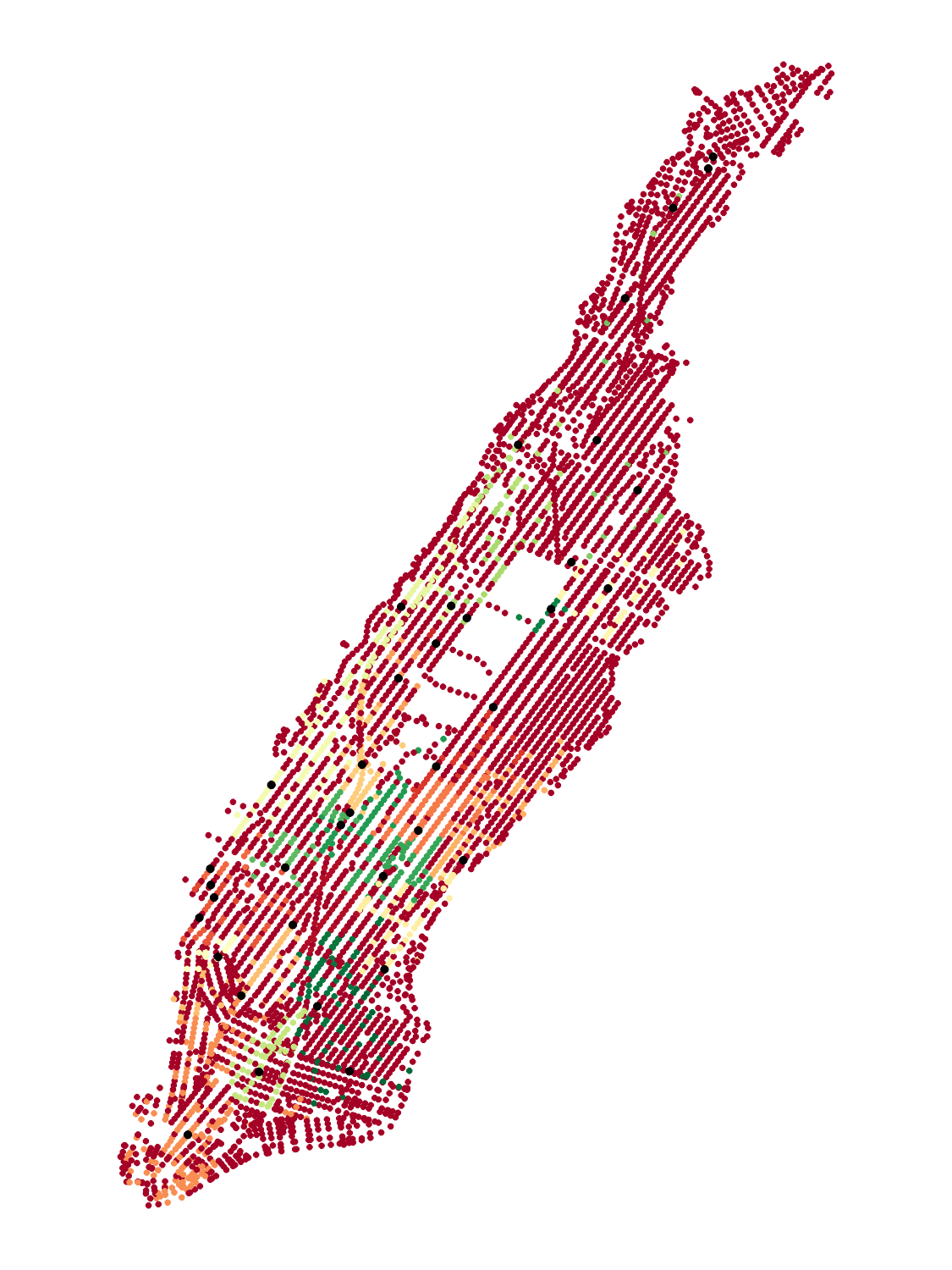}}
  \hfill
  \vspace{-10pt}
  \caption{Partition with different methods}

\label{fig:partition_algorithms}
\end{figure}

\textbf{Partition with Different Algorithms:} Fig. \ref{fig:partition_algorithms} gives the partitioning results using different methods. Here the different color of vertices shows that it belongs to different sub-graphs. The black vertices stand for the relocation centers. In Fig. \ref{fig:partition_wkm}, we use ReLU activation function in the objective function. It is clear that partitioning with four baseline methods, we get the sub-graphs with clear boundaries and the size of each sub-graph does not have significant differences. However, with the dFDA-VeD, a huge amount of vertices belongs to the same sub-graphs. This is reasonable as many vertices are no demand or low demand points (as shown in Fig. \ref{fig:demands_gap}). This give us the insight that the graph partition for relocation searching do not need to be similar size or have clear boundary.  

\begin{figure}[!ht]

\centering
  \subfigure[ignore]{
  \label{fig:partition_ignore}
  \includegraphics[width=0.15\linewidth]{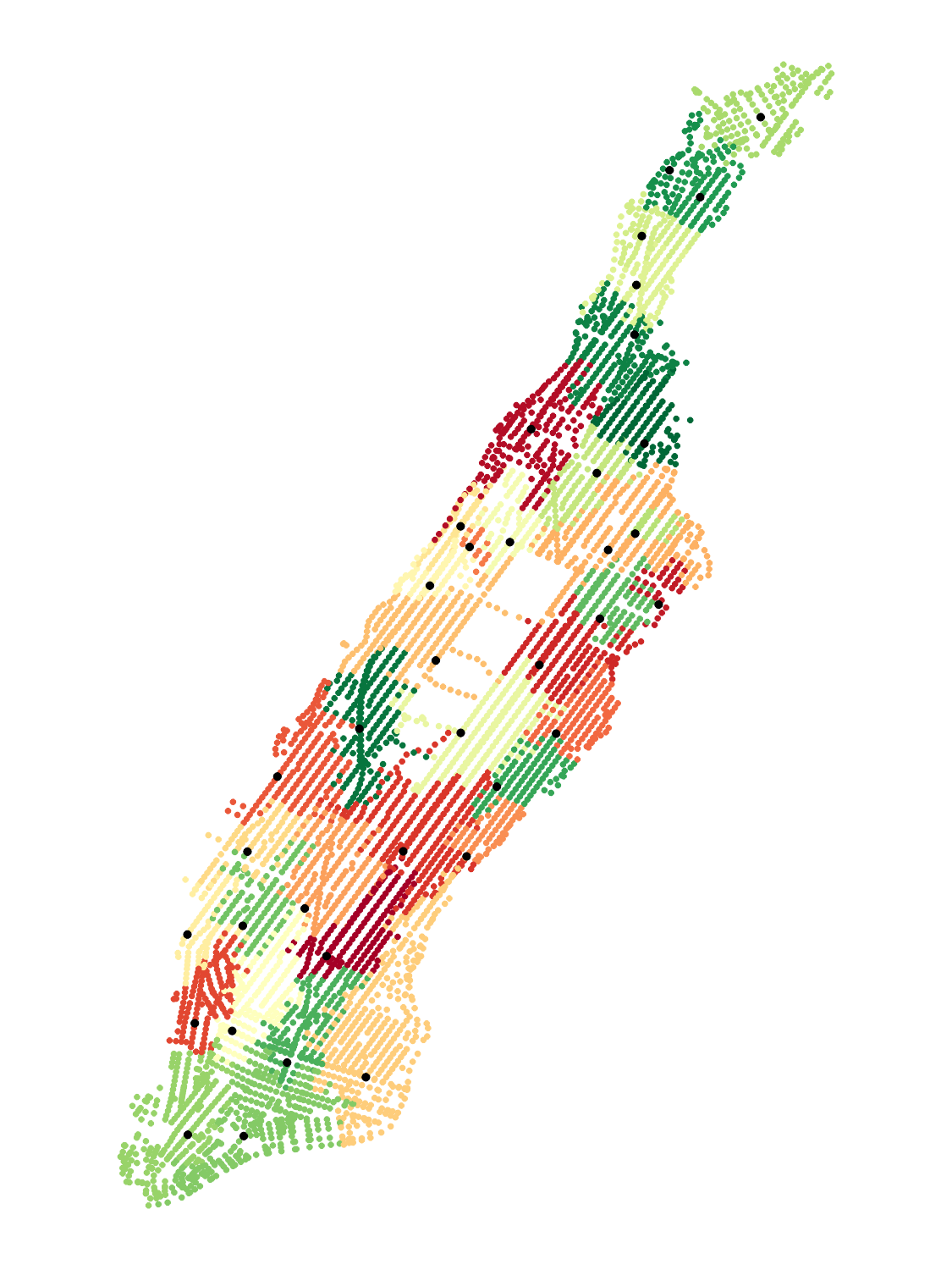}}
  \hfill
  \subfigure[identity]{
  \label{fig:partition_identity}
  \includegraphics[width=0.15\linewidth]{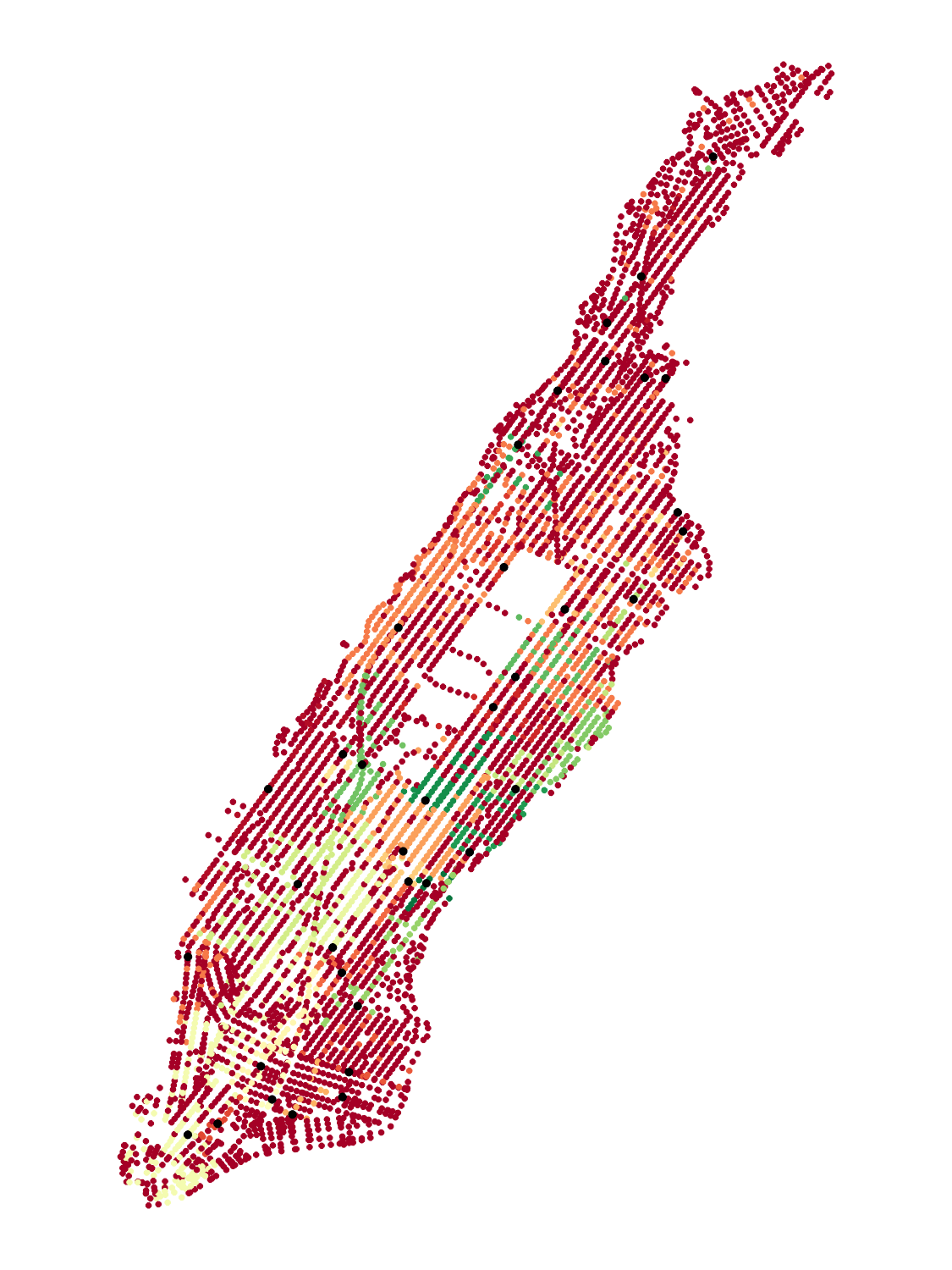}}
  \hfill
  \subfigure[Sigmoid]{
  \label{fig:partition_sigmoid}
  \includegraphics[width=0.15\linewidth]{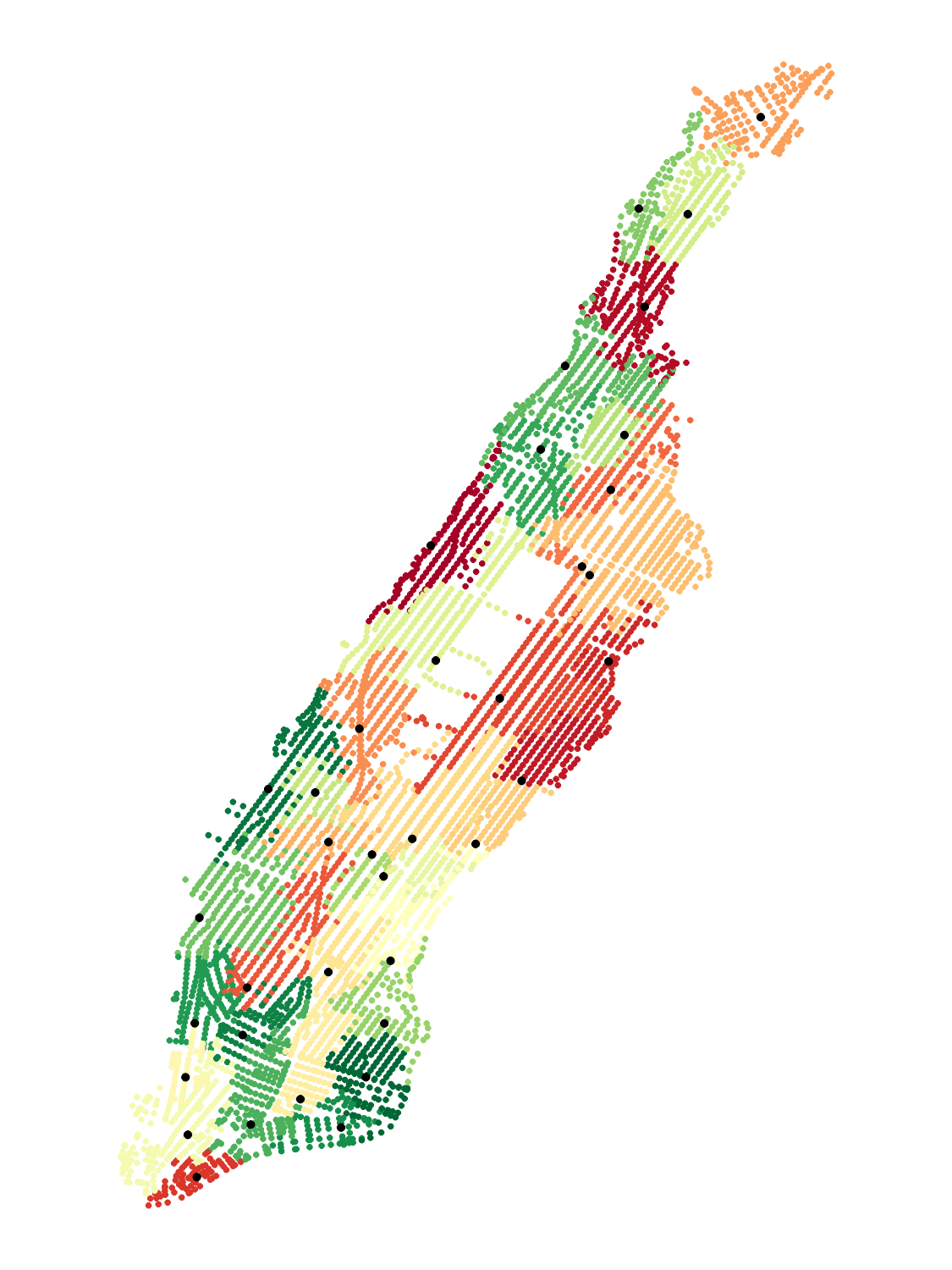}}
  \hfill
  \subfigure[Softplus]{
  \label{fig:partition_softplus}
  \includegraphics[width=0.15\linewidth]{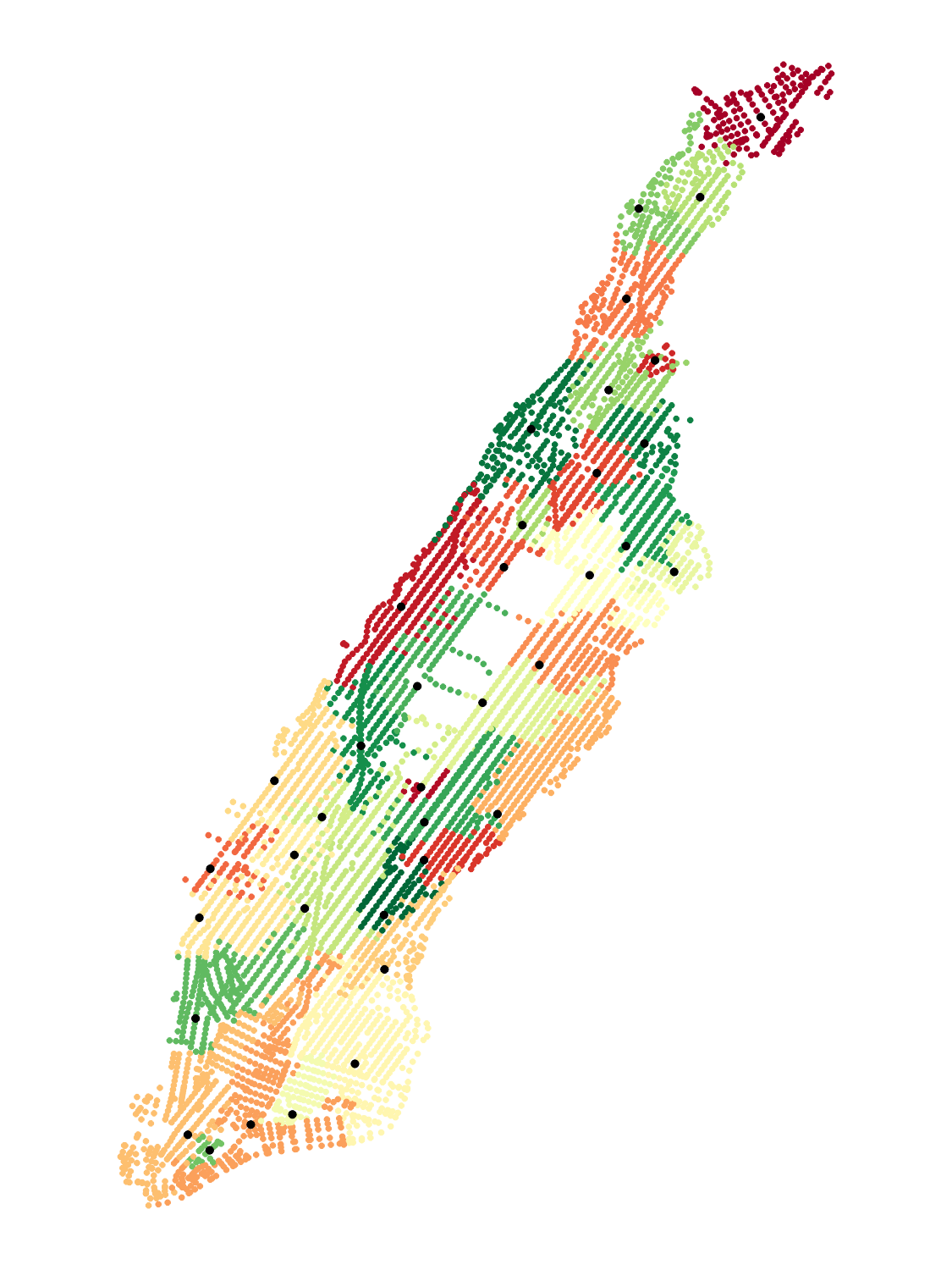}}
  \hfill
  \subfigure[ReLU]{
  \label{fig:partition_ReLU}
  \includegraphics[width=0.15\linewidth]{figures/partitionReLU-eps-converted-to.pdf}}
  \hfill
  \vspace{-10pt}
  \caption{Partition with different activation functions}
  
\label{fig:partition_activations}
\end{figure}

\textbf{Partition with Different Activation Functions:} Fig \ref{fig:partition_activations} gives the partitioning results using relocation center searching algorithm with different activation functions. It shows that with the Identity and ReLU function similar partitioning result can be achieved. This is consistent with their similar performance in terms of served trip ratio $R$ (as shown in Fig. \ref{fig:ratios_activation}). This give us the insight that we should not treat the high demand and low demand as the same importance (just like Ignore function does), and should not treat the low demand points in a positive way (just like Sigmoid and Softplus function do). The low demand and high demand area should be treated differently. With the proper way (e.g. ReLU and Identity function) to treat different demand, the idle vehicle relocation could significantly improve the served trip ratio.

\section{Conclusion}
\label{sec:conclusion}

In this paper, we proposed a dynamic future demand aware vehicle dispatching system, called \texttt{dFDA-VeD}. The proposed system is based on relocating the idle vehicles to re-balance the sub-areas in an urban region. As the traffic conditions and travel demands continuously change in a dynamic manner, the relocation is done by dynamically identifying the relocation centres taking the real-time conditions into account. We demonstrate the performance of our dynamic future demand aware vehicle dispatching system through extensive experiments on real data. We outperform the existing state-of-the-art methods and vehicle dispatching systems in terms of serving ratio. An important future research direction is to develop an index for an efficient maintenance of the dynamic information to support decision making in the dispatching system.



\bibliographystyle{ACM-Reference-Format}
\bibliography{reference}
%
\end{document}